\pdfoutput=1

\documentclass[12pt,a4paper]{article}
\usepackage[margin=1in]{geometry}
\usepackage[sort&compress,square,numbers]{natbib}
\usepackage{lipsum,times,graphicx}
\usepackage[labelfont=bf]{caption}
\usepackage[T1]{fontenc}
\usepackage{array}
\usepackage{multirow}
\usepackage{amsmath}
\usepackage{amssymb}
\usepackage{stmaryrd}
\usepackage{graphicx}
\usepackage{wasysym}
\usepackage{diagbox}
\usepackage{lipsum}
\usepackage{amsthm}
\usepackage{enumerate}
\usepackage[shortlabels]{enumitem}
\usepackage[colorlinks,
linkcolor=blue,
anchorcolor=blue,
citecolor=blue]{hyperref}
\usepackage{float}
\usepackage{graphicx}
\usepackage{placeins}
\usepackage{makecell}
\usepackage{xr}
\usepackage{cleveref}

\theoremstyle{definition}
\newtheorem{thm}{Theorem}
\newtheorem{defn}{Definition}

\newtheorem{prop}{Proposition}
\newtheorem{rem}{Remark}

\newtheorem{conj}{Conjecture}

\usepackage[english]{babel}

\numberwithin{equation}{section}

\begin{document}

\title{On rigid origami \uppercase\expandafter{\romannumeral3}: local rigidity analysis}

\author{
	Zeyuan He$^{1}$, Simon D. Guest$^{2}$}

\author{Zeyuan He, Simon D. Guest\footnote{zh299@cam.ac.uk, sdg@eng.cam.ac.uk. Civil Engineering Building, Department of Engineering, 7a JJ Tomson Ave, University of Cambridge, Cambridge CB3 0FA, United Kingdom}}
\maketitle

\begin{abstract}
\noindent Rigid origami is examined from the perspective of rigidity theory.  First and second order rigidity are defined from local differential analysis of the consistency constraint; while the static rigidity and prestress stability are defined after finding the form of internal force and load. This article will show that first-order or static rigidity implies prestress stability, which implies second-order rigidity, which implies rigidity, but none of these is reversible.  Examples are given of rigid origami structures with these different kinds of rigidity.  Examining the different aspects of the rigidity of origami might give a novel perspective for the development of new folding patterns, or for the design of origami structures where some rigidity is required.
\end{abstract}

\begin{small}
	{$\; \; \,$ \textbf{Keywords}:} foldability, stress, load, first-order, prestress stability, second-order
\end{small}


	
\section{Introduction}

Rigid origami has been developed as a tool for effectively transforming a two-dimensional material into a three-dimensional structure, hence most of previous studies focus on the kinematics and mechanical properties of foldable rigid origami. In this article, we will consider a different viewpoint -- a rigid origami that is not foldable, among which there is a hierarchical relation of different levels of rigidity. The local rigidity concepts, including the first-order or static rigidity, prestress stability and second-order rigidity are similar to those used for classical bar-joint frameworks, but there are also some special features. Considering the different aspects of the rigidity for rigid origami might give a novel perspective for the development of new folding patterns, or for the design of origami structures where some rigidity is required.

\begin{figure}[hbtp]
	\noindent \begin{centering}
		\includegraphics[width=1\linewidth]{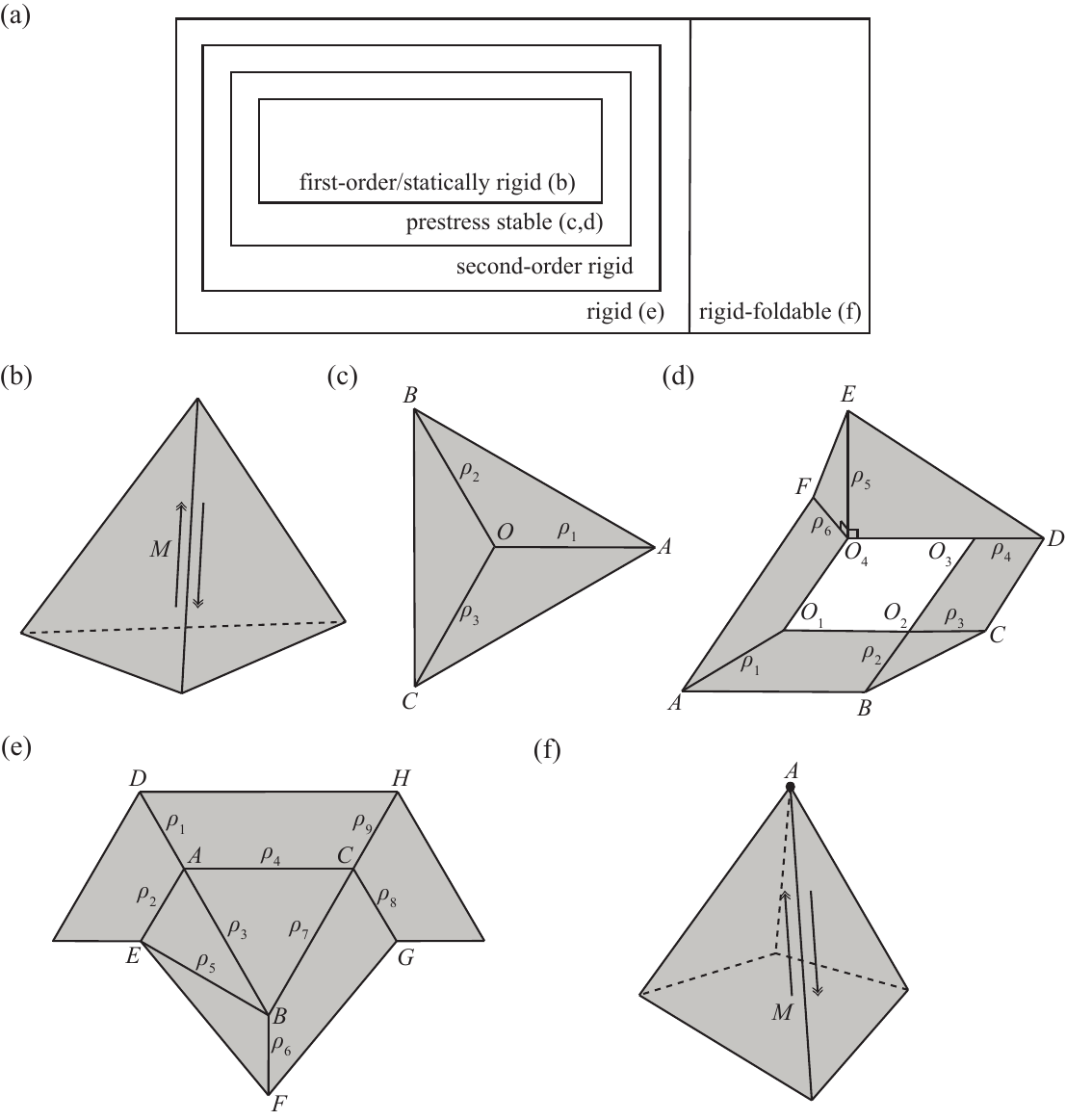}
		\par\end{centering}
	
	\caption{\label{fig: introduction} (a) shows the hierarchical relation among different levels of rigidity for a rigid origami. The first-order rigidity is equivalent to static rigidity, which implies prestress stability, which implies second-order rigidity, which implies rigidity, but none of these relationships is reversible. (b)--(f) are examples corresponding to each region in (a)  --- vertices are shown by capital letters, and folding angles by $\rho_{i}$. (b) is a tetrahedron that is first-order (equivalently, statically) rigid. (c) is a planar degree-3 vertex, (d) is a non-planar degree-6 hole, both of which are prestress stable but not first-order rigid. (e) is a planar 3-vertex rigid origami that is rigid but not second-order rigid. (f) is a rigid-foldable degree-4 vertex. $M$ is a pair of opposite torques applied on two panels incident to an inner crease, which is the form of load for rigid origami. (b) is able to carry the load in this configuration, while (f) cannot.}
\end{figure}

The relationship between rigidity concepts, as well as examples for each level of rigidity, are presented in figure~\ref{fig: introduction}.  The hierarchical relation in (a) will be proved in Section \ref{section: relation}. The concepts underlying the examples in (b)--(f) will be presented in the rest of the paper. The tetrahedron shown in (b) is first-order rigid (described in Section~\ref{section: first-order}), i.e., the only first-order flex in the tangent space of the consistency constraint is $\boldsymbol{0}$.  The tetrahedron is also statically rigid (described in Section~\ref{section: statics}), i.e., under any load applied, there would be a set of internal forces generated to keep the rigid origami in equilibrium. Here, the form of load is a pair of opposite torques applied on two panels incident to each inner crease. The example in (c) is the simplest rigid origami to be prestress stable (described in Section~\ref{section: prestress stability}) but not first-order rigid. There is a one-dimensional first-order flex $\{\rho_1',\rho_2',\rho_3'\}=a_1\{1,1,1\}$, $a_1 \in \mathbb{R}$. Although this vertex would be "shaky" along this first-order flex, this configuration will reach a strict local minimum of a predefined potential energy function with a positive self-stress, therefore it will not be able to deform greatly. Example in (d) has two-dimensional first-order flex and self-stress, which is also examined to be prestress stable. It turns out that prestress stability is a relatively strong class for common rigid origami structures that are not first-order rigid. The next level is second-order rigidity (described in Section~\ref{section: second order rigidity}). There is no self-stress that could help the rigid origami reach a local minimum of potential energy, and there is also no first-order flex that can be extended to a second-order flex. Second-order rigidity will imply rigidity. The example in (e) is rigid but not second-order rigid, which relies on special choice of sector angles. The example in (f) is a degree-4 single-vertex with centre vertex $A$ where the sum of sector angles is less than $2\pi$, and this vertex is rigid-foldable. Only some loads can be carried. Here we show a load that acts to change the configuration of this single-vertex. An interesting point is for some special rigid origami, different levels of rigidity might be equivalent. Such extension of local rigid-foldability is not easy to predict, and this will be discussed in Section~\ref{section: extension}.

The inspiration for the levels of rigidity comes from classical studies on the statics, prestress stability and second-order rigidity of bar-joint frameworks, e.g.\ \cite{connelly_second-order_1996,connelly_frameworks_2015,schulze_rigidity_2017} --- the preliminaries for these materials are provided in Section S1 of the supplementary material. We find a good correspondence between rigid origami and bar-joint framework. However, rigid origami has some special features. Because of our kinematic definitions, we don't need to consider the Euclidean motion of a rigid origami in the folding angle expression - the only trivial flex is $\boldsymbol{0}$, which simplifies some conclusions. Second, in line with classical rigidity theory for a framework, a rigid origami has its special form of the underlying graph, body and linkage in the folding angle expression. For a bar-joint framework, the Jacobian and Hessian of consistency constraints are linear and constant, while for a rigid origami they are in a totally different form of higher order. As a consequence, the form of internal force and load are also different. Third, the effect of self-intersection could not be revealed from the classic method of doing an algebraic analysis on the consistency constraints. The collision between panels might induce rigidity, but is not considered here --- numeric methods are more likely to be efficient when dealing with self-intersection. Fourth, a set of folding angles might correspond to several stacking sequences (an example is in \cite[figure 3]{he_rigid_2019}), and different stacking sequences might also behave differently when considering the self-intersection of panels, which we also do not consider in this paper. Fifth, when some folding angles are $\pm \pi$, a flex is valid only when it points away from $\pm \pi$ since the range of a folding angle is $[-\pi,\pi]$. Some examples on this topic are provided in Section S3.6 of the supplementary material. In this article we will focus on the local algebraic analysis of rigid origami. The effect of self-intersection and stacking sequence, and the ``one-side'' property when some folding angles are $\pm \pi$, are some topics that require further work. In the rest of this paper we will require every folding angle to be in $(-\pi,\pi)$, but it does not mean the conclusions drawn for local rigidity will always fail when some folding angles are $\pm \pi$.

This paper considers rigidity to second order.  It is natural to ask ``could this hierarchical relation be extended to countable order of rigidity and ends at finite rigidity?''. However, there seems to exist a limit for such local differential analysis. It turns out that there exists a bar-joint framework that is third-order rigid and flexible \cite{connelly_higher-order_1994}, which implies the chain relation of local rigidity may not be closed by rigidity, or we might need to modify the definition of local rigidity. Further, for a bar-joint framework, a sufficiently high order flexibility will be equivalent to flexibility; or with some extra conditions, they would be equivalent even if the order is not that high \cite{alexandrov_sufficient_1998}. There is a proposal for a revised definition of local rigidity \cite{stachel_proposal_2007}, but a complete theory still requires certain amount of work.

\section{Consistency constraint on folding angles}

In this section we will briefly recap some basic definitions on rigid origami and the algebraic constraint on folding angles. The detailed version is provided in \cite{he_rigid_2019}.

\begin{figure}[!tb]
\noindent \begin{centering}
	\includegraphics[width=1\linewidth]{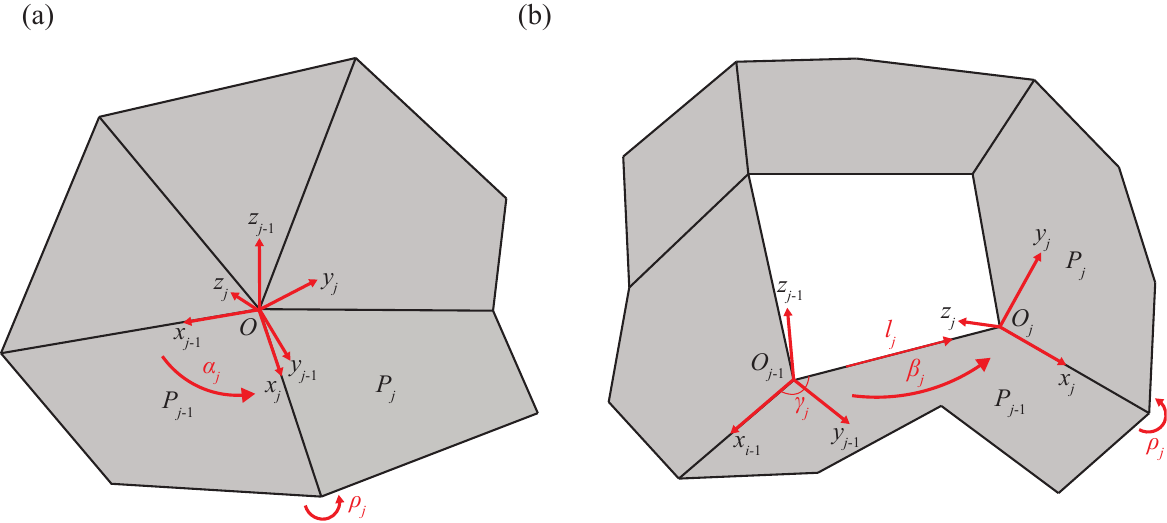}
	\par\end{centering}

\caption{\label{fig: single vertex and hole}(a) and (b) show the rotation and possible translation of local coordinate systems around a degree-5 vertex and a degree-5 hole. The first Betti number of (a) is 0 and of (b) is 1. The form of consistency constraint is with equation \eqref{eq: consistency 1} for (a) and with equation \eqref{eq: consistency 2} for (b).}
\end{figure}

\begin{defn}
An \textit{underlying graph} $G$ is a \textit{multi-level graph} where an edge connects multiple vertices in sequence. A \textit{realization} $G(\boldsymbol{\rho},\lambda')$ of the underlying graph $G$ is a panel-hinge framework where each hinge (also called an inner crease) corresponds to a vertex in $G$, and inner creases connected by an edge in $G$ form an inner vertex or hole of the panel-hinge framework in the sequence described by this edge in $G$. Here $G$ can be regarded as a dual of the interior of the crease pattern. Each inner crease is the intersection of two panels, and a panel is a polygon bounded by creases. The \textit{consistency constraint} around each inner vertex or hole is a series of rotations and translations (will be explained below) represented by the sector angles and folding angles. A \textit{sector angle} $\alpha \in (0, 2\pi) $ is the angle between adjacent inner creases on a panel. A \textit{folding angle} $\rho \in [-\pi, \pi]$ on an inner crease is the angle measuring how its adjacent two panels deviate from a plane viewed from a given orientation. The collection of sector angles and folding angles for all inner creases are written as $\boldsymbol{\alpha}$ and $\boldsymbol{\rho}$. For a rigid origami $G(\boldsymbol{\rho})$, there might be multiple stacking sequences corresponding to the same folding angle $\boldsymbol{\rho}$, therefore we introduce the \textit{order function} $\lambda'$ defined on possible non-crease contact points on non-adjacent panels to describe the stacking sequence. The order function $\lambda'$ must satisfy several conditions to prevent self-intersection.   

Given a rigid origami $G(\boldsymbol{\rho},\lambda')$, a \textit{rigidly folded state} is another realization of $G$ with the same sector angles $\boldsymbol{\alpha}$, which is also an isometry $f$ excluding Euclidean motion of the given rigid origami. A \textit{rigid folding motion} is a family of continuous isometry mapping each time $t \in [0,1]$ to a rigidly folded state, where the order function $\lambda'$ should also satisfy the continuity condition to guarantee the rigid folding motion to be physically admissible without self-intersection. If there is a rigid folding motion starting from a rigid origami, it is \textit{rigid-foldable}, otherwise \textit{rigid}.
\end{defn}
Next we will clarify what is the constraint on folding angles $\boldsymbol{\rho}$ if regarding the sector angles $\boldsymbol{\alpha}$ as known parameters. The order function $\lambda'$ is only valid when there is contact between panels not adjacent to an inner crease, and could be examined later. The algebraic constraint on folding angles is the collection of consistency condition on each inner vertex and hole, which is derived from the rotation and translation of local coordinate systems on each panel near a vertex or a hole (figure \ref{fig: single vertex and hole}). This is the sufficient and necessary condition for a set of folding angles to be an element of the configuration space if allowing the self-intersection of panels. Here we further clarify that a \textit{hole} refers to a \textit{homology class} in the \textit{first homology group}, and the number of holes is called the \textit{first Betti number} \cite{hatcher_algebraic_2005}. Homology itself was developed as a way to analyse and classify manifolds according to their cycles -- closed loops that can be drawn on a given manifold but can not be continuously deformed into each other. Informally, cycles that can be continuously transformed into each other belong to the same homology class of the first homology group. A hole in the consistency condition is an arbitrary representative cycle of a homology class on the paper. The first Betti number is also the maximum number of cuts that can be made without dividing a surface into two separate pieces. For example, the first Betti number of a sphere and a disk is 0; of a cylindrical surface is 1; of a torus is 2. 

The local coordinate systems are built in the following way when deriving the consistency constraint. (a) Around a vertex surrounded by $n$ panels (figure \ref{fig: single vertex and hole}(a)), a \textit{local coordinate system} is built on each panel $j$ ($j \in [1,n]$), whose origin $O_j$ is on the centre vertex, $x$-axis is on an inner crease, pointing outside the origin, $z$-axis is normal to the panel. The direction of all $z$-axes of local coordinate systems are consistent with the orientation of the paper and hence consistent with the definition of the sign of folding angles. Specifically, the transformation between local coordinate systems of panel $j-1$ and panel $j$ is a rotation $\alpha_j$ along $z_{j-1}$, and a rotation $\rho_j$ along $x_j$. After a series of rotations the coordinate system returns to the one built on panel $n$. The matrix form of transformation is given in equation \eqref{eq: consistency 1}. (b) Around a hole surrounded by $n$ panels (figure \ref{fig: single vertex and hole}(b)), we build the local coordinate systems similarly. Each origin $O_j$ is on a boundary vertex of this hole, $x$-axis is on an inner crease, pointing outside each origin, $z$-axis is normal to the panel. The transformation between local coordinate systems of panel $j-1$ and panel $j$ is a translation $[l_j \cos \gamma_j; l_j \sin \gamma_j; 0]$ measured in the coordinate system built on panel $j-1$, followed by a rotation $\beta_j$ along $z_{j-1}$, and a rotation $\rho_j$ along $x_j$. $\beta_j$ and $\gamma_j$ can be linearly expressed by the sector angles $\boldsymbol{\alpha}$. The matrix form is given in equation \eqref{eq: consistency 2}.

At every inner vertex $v_i~(1 \le i \le N_v$, figure \ref{fig: single vertex and hole}(a)):
\begin{align} \label{eq: consistency 1}
R(\boldsymbol{\rho})=\prod_{1}^{\mathrm{deg}(v_i)}
\left[ \begin{array}{ccc}
	\cos \alpha_{j} & -\sin \alpha_{j} & 0\\
	\sin \alpha_{j} &  \cos \alpha_{j} & 0\\
	0                &                 0 & 1
\end{array} \right]
\left[ \begin{array}{ccc}
	1                &                 0 & 0             \\
	0                &     \cos \rho_{j} & -\sin \rho_{j}\\
	0                &     \sin \rho_{j} &  \cos \rho_{j}
\end{array} \right]=I
\end{align}
where $N_v$ is the number of inner vertices, $\mathrm{deg}(v_i)$ is the number of creases (degree) incident to $v_i$, $\alpha_{j}$ is between axes $x_{j-1}$ and $x_j$ ($2 \le j \le \mathrm{deg}(v_i)$), $\alpha_1$ is between axes $x_{\mathrm{deg}(v_i)}$ and $x_1$. $R$ is formed by post-multiplication. Only three of the nine equations are independent, which are in different columns and rows. 

At every hole with boundary $h_i~(1 \le i \le N_h$, figure \ref{fig: single vertex and hole}(b)):
\begin{align} \label{eq: consistency 2}
T(\boldsymbol{\rho})=\prod_{1}^{\mathrm{deg}(h_i)}
\left[ \begin{array}{cccc}
	\cos \beta_j & -\sin \beta_j & 0 & l_j \cos \gamma_j \\
	\sin \beta_j & \cos \beta_j & 0 & l_j \sin \gamma_j \\
	0 & 0 & 1 & 0   \\
	0 & 0 & 0 & 1
\end{array} \right]
\left[ \begin{array}{cccc}
	1 & 0 & 0 & 0 \\
	0 & \cos \rho_j & -\sin \rho_j & 0 \\
	0 & \sin \rho_j & \cos \rho_j & 0 \\
	0 & 0 & 0 & 1
\end{array} \right]=I
\end{align}
where $N_h$ is the number of holes, $\mathrm{deg}(h_i)$ is the number of creases (degree) incident to $h_i$, $\beta_{j}$ is between axes $x_{j-1}$ and $x_j$ ($2 \le j \le \mathrm{deg}(h_i)$), $\beta_1$ is between axes $x_{\mathrm{deg}(h_i)}$ and $x_1$. $[l_j \cos \gamma_j,l_j \sin \gamma_j,0]$ ($1 \le j \le \mathrm{deg}(h_i)$) is the position of $O_j$ measured in the local coordinate system for panel $j-1$. $T$ is formed by post-multiplication. Only six of the sixteen equations are independent. Three of them are in the top left $3 \times 3$ rotation matrix, the other three are the elements from row 1 to row 3 in column 4. If the inner creases are concurrent, this hole will degenerate to a vertex. 

If there are $N_v$ inner vertices and $N_h$ holes, the number of independent consistency constraint will be $3N_v+6N_h$. It will simplify further algebraic analysis if considering the independent components of the consistency constraint. Why choosing the particular independent components below will be explained in Section \ref{section: first-order}.  

For each vertex:
\begin{align} \label{eq: independent consistency 1}
\left[ \begin{array}{ccc}
	\ast & \ast & A_2 \\
	A_3 & \ast & \ast \\
	\ast & A_1  & \ast
\end{array} \right]=R(\boldsymbol{\rho})
\end{align}
For each hole:
\begin{align} \label{eq: independent consistency 2}
\left[ \begin{array}{cccc}
	\ast & \ast & A_2 & A_4 \\
	A_3 & \ast & \ast & A_5 \\
	\ast & A_1  & \ast & A_6   \\
	0 & 0 & 0 & 1
\end{array} \right]=T(\boldsymbol{\rho})
\end{align}
Here $\ast$ means elements that are not important in further discussion. $\boldsymbol{A}$ is a vector with length $3N_v+6N_h$ assembled from 3 components for each vertex and 6 components for each hole. 

If $\boldsymbol{\rho}$ is a solution of consistency constraint, $\boldsymbol{A}(\boldsymbol{\rho})=\boldsymbol{0}$. However, the converse is not necessarily true: if $\boldsymbol{A}(\boldsymbol{\rho})=\boldsymbol{0}$, $\boldsymbol{\rho}$ might not be a solution of the consistency constraint, because equations \eqref{eq: independent consistency 1} and \eqref{eq: independent consistency 2} could give some rotation matrices whose determinant is 1 but formed by 0 and $\pm 1$ apart from the Identity. In other words, the solution space of the independent components of the consistency constraint $\boldsymbol{A}(\boldsymbol{\rho})=\boldsymbol{0}$ is larger than the solution space of the consistency constraint. However, these solutions can be easily removed by examination. The first and second order derivative of the independent consistency constraint will be used in the analysis of first-order rigidity, prestress stability and second-order rigidity in Sections \ref{section: first-order}, \ref{section: prestress stability} and \ref{section: second order rigidity}.

\section{First-order rigidity} \label{section: first-order}

The first-order rigidity and a first-order flex are defined as below. 

\begin{defn}
A rigid origami $(\boldsymbol{\rho}, \lambda')$ is \textit{first-order rigid} if the only solution of $\mathrm{d} \boldsymbol{A}/\mathrm{d} \boldsymbol{\rho} \cdot \boldsymbol{\rho}'= \boldsymbol{0}$ with respect to $\boldsymbol{\rho}'$ is $\boldsymbol{0}$, equivalently, the rank of \textit{rigidity matrix} $\mathrm{d} \boldsymbol{A}/\mathrm{d} \boldsymbol{\rho}$ equals to the number of inner creases $N_c$. Otherwise this rigid origami is \textit{first-order rigid-foldable}. A non-zero $\boldsymbol{\rho}'$ is called a \textit{first-order flex}, which forms a linear space of dimension $N_c-\mathrm{rank}(\mathrm{d} \boldsymbol{A}/\mathrm{d} \boldsymbol{\rho})$. 
\end{defn}

We will show how to derive the rigidity matrix $\mathrm{d} \boldsymbol{A}/\mathrm{d} \boldsymbol{\rho}$ for a large rigid origami after writing $\mathrm{d} \boldsymbol{A}/\mathrm{d} \boldsymbol{\rho}$ for its restriction on a single-vertex or single-hole.

\subsection{Rigidity matrix for a single-vertex or hole}

Consider the first-order derivative of equations \eqref{eq: independent consistency 1} and \eqref{eq: independent consistency 2}.
\begin{equation} \label{eq: Jacobian 1}
\dfrac{\partial }{\partial \rho_j} 	\left[ \begin{array}{ccc}
	\ast & \ast & A_2 \\
	A_3 & \ast & \ast \\
	\ast & A_1  & \ast
\end{array} \right]= \left[ \begin{array}{cccc}
	0 & -x_{3j} & x_{2j} \\
	x_{3j} & 0 & -x_{1j} \\
	-x_{2j} & x_{1j} & 0 
\end{array} \right]
\end{equation}
where $\boldsymbol{x}_j=[x_{1j};x_{2j};x_{3j}]$ is the direction (column) vector of the inner crease $\rho_j$ measured in a global coordinate system, pointing away from this vertex.
\begin{equation} \label{eq: Jacobian 2}
\dfrac{\partial }{\partial \rho_j} \left[ \begin{array}{cccc}
	\ast & \ast & A_2 & A_4 \\
	A_3 & \ast & \ast & A_5 \\
	\ast & A_1  & \ast & A_6   \\
	0 & 0 & 0 & 1
\end{array} \right] = 
\left[ \begin{array}{cccc}
	0 & -x_{3j} & x_{2j} & \multirow{3}{*}{$\boldsymbol{O}_j \times \boldsymbol{x}_j$}\\
	x_{3j} & 0 & -x_{1j} & \\
	-x_{2j} & x_{1j} & 0 & \\
	0 & 0 & 0 & 1
\end{array} \right]
\end{equation}
where $\boldsymbol{x}_j=[x_{1j};x_{2j};x_{3j}]$ is the direction (column) vector of the inner crease $\rho_j$ measured in a global coordinate system, pointing away from this hole. $\boldsymbol{O}_j$ is the position of vertex on the hole incident to $\rho_j$ measured in the global coordinate system (figure \ref{fig: single vertex and hole}). The derivation is provided in Section S2 of the supplementary material. 

The reason for choosing the particular components when defining the independent consistency constraint $\boldsymbol{A}(\boldsymbol{\rho})=\boldsymbol{0}$ is to ensure that we capture independent non-zero values in equations \eqref{eq: Jacobian 1} and \eqref{eq: Jacobian 2} to describe the ``speed'' of a dynamic system. 

The matrix form of rigidity matrix for a degree-$n$ single-vertex or single-hole could therefore be written as:
\begin{equation} 	\label{eq: Jacvertex}
\dfrac{\mathrm{d} \boldsymbol{A}}{\mathrm{d} \boldsymbol{\rho}}_\text{vertex}=
\left[{\begin{array}{cccc}
		\boldsymbol{x}_1 & \boldsymbol{x}_2 & \cdots & \boldsymbol{x}_n
\end{array}} \right],
\end{equation}
\begin{equation} 	\label{eq: Jachole}
\dfrac{\mathrm{d} \boldsymbol{A}}{\mathrm{d} \boldsymbol{\rho}}_\text{hole}=
\left[{\begin{array}{cccc}
		\boldsymbol{x}_1 & \boldsymbol{x}_2 & \cdots & \boldsymbol{x}_{n}\\
		\boldsymbol{O}_1 \times \boldsymbol{x}_1 & \boldsymbol{O}_2 \times \boldsymbol{x}_2  & \cdots & \boldsymbol{O}_n \times \boldsymbol{x}_n
\end{array}} \right]
\end{equation}

The rigidity matrix $\mathrm{d} \boldsymbol{A}/\mathrm{d} \boldsymbol{\rho}$ for each single vertex or hole can also be explained from analytical mechanics. First, around each degree-$n$ single vertex, the virtual rotation from panel $n$ to panel $j$ induced by a perturbation on folding angles $\delta \boldsymbol{\rho}$ is $\boldsymbol{x}_1 \delta \rho_1+\boldsymbol{x}_2 \delta \rho_2+...+\boldsymbol{x}_j \delta \rho_j$. After returning to panel $n$, the relative virtual rotation should be $0$, therefore,
\begin{equation} \label{eq: virtual 1}
\sum_{1}^{n} \boldsymbol{x}_j \delta \rho_j=\boldsymbol{0}
\end{equation}
Around each degree-$n$ hole, equation \eqref{eq: virtual 1} still holds. Fix panel $n$ to exclude Euclidean motion, the virtual displacement from origin of the global coordinate system $\boldsymbol{0}$ to the $xy$-plane of local coordinate system built on panel $j$ induced by a perturbation on folding angles $\delta \boldsymbol{\rho}$ is $\boldsymbol{O}_1 \times \boldsymbol{x}_1 \delta \rho_1 + \boldsymbol{O}_2 \times \boldsymbol{x}_2 \delta \rho_2 +...+ \boldsymbol{O}_j \times \boldsymbol{x}_j \delta \rho_j$. For panel $n$, this virtual displacement should be $0$, hence around each hole we have
\begin{equation} \label{eq: virtual 2}
\sum_{1}^{n} (\boldsymbol{O}_j \times \boldsymbol{x}_j) \delta \rho_j  =\boldsymbol{0}
\end{equation}

\subsection{Measure of deformation about a vertex or hole} \label{subsection: measure}

\begin{figure}
\noindent \begin{centering}
	\includegraphics[scale=0.6]{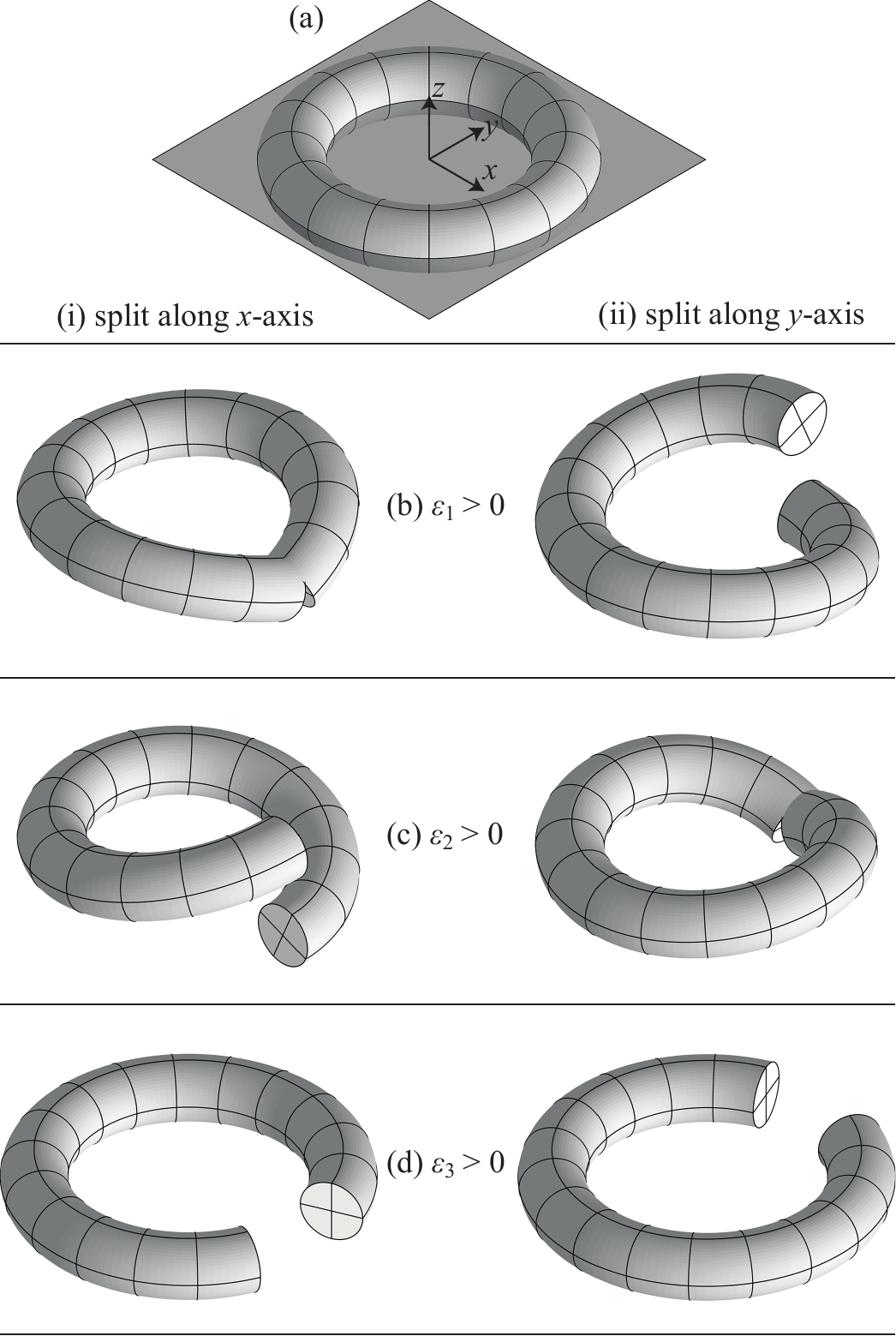}
	\par\end{centering}
\caption{\label{fig: ring}An illustration of the first-order error $\boldsymbol{\varepsilon}$ in a global coordinate system around a single vertex. Note that there is a fundamental difficulty in showing this deformation, as it is a measure of what rigid origami \emph{cannot} do, so a cut has to be introduced to allow the deformation to be shown, and where the cut is introduced \textit{may} make a big difference to the appearance. (a) Using a flat vertex as an example, we attach a torus to the paper around the vertex.  (b) We first consider a deformation $\varepsilon_1$ distributed evenly around the torus, so that there is a constant `curvature' around the $x$-axis (which in places manifests as a twist). The torus is then not able to close, so we show the deformed torus with a cut along the (i) $x$-axis, or (ii) $y$-axis. Although the images look very different, the underlying deformation is the same in each case.  (c,d) We similarly show the deformation $\varepsilon_2$ and $\varepsilon_3$.  Only for (b)(i) and (c)(ii) would a first-order error be compatible with the rotation of a hinge, in $x$- and $y$-directions, respectively.  The internal forces that are work-conjugate with these first-order errors are shown in figure \ref{fig: intforces}.}
\end{figure}

Given a rigid origami $(\boldsymbol{\rho},\lambda')$, it is instructive to consider the form of the deformations that we are not allowing if there is a perturbation on folding angles $\delta \boldsymbol{\rho}$. A measure of deformation called \textit{first-order error} $\boldsymbol{\varepsilon}(\boldsymbol{\rho})$ could be derived from the first-order estimation of the independent consistency constraint $\boldsymbol{A}(\boldsymbol{\rho}+\delta \boldsymbol{\rho})$:
\begin{equation} \label{eq: first-order deformation}
\boldsymbol{\varepsilon}(\boldsymbol{\rho})=\dfrac{\mathrm{d} \boldsymbol{A}}{\mathrm{d} \boldsymbol{\rho}} \delta \boldsymbol{\rho}
\end{equation}
For a degree-$n$ vertex:
\begin{equation} \label{eq: error 3}
\begin{aligned}
	\left[\begin{array}{c}
		A_1 \\
		A_2 \\
		A_3
	\end{array}\right] & = \left[\begin{array}{c}
		\varepsilon_1 \\
		\varepsilon_2 \\
		\varepsilon_3
	\end{array}\right] +O(\delta \boldsymbol{\rho}^2) =\sum_{1}^{n} \boldsymbol{x}_j \delta \rho_j+O(\delta \boldsymbol{\rho}^2)
\end{aligned}
\end{equation}
Here the first-order error is the components about the global $x$, $y$ and $z$ axes of the rotation from the local coordinate system built on panel $n$ to itself, as a circuit is taken around the vertex with folding angles $\boldsymbol{\rho}+\delta \boldsymbol{\rho}$.

For a degree-$n$ hole:
\begin{equation} \label{eq: error 4}
\begin{aligned}
	\left[\begin{array}{c}
		A_1 \\
		A_2 \\
		A_3 \\
		A_4 \\
		A_5 \\
		A_6
	\end{array}\right] & =\left[\begin{array}{c}
		\varepsilon_1 \\
		\varepsilon_2 \\
		\varepsilon_3 \\
		\varepsilon_4 \\
		\varepsilon_5 \\
		\varepsilon_6
	\end{array}\right]+O(\delta \boldsymbol{\rho}^2) =\left[\begin{array}{c}
		\sum \limits_{1}^{n} \boldsymbol{x}_j \delta \rho_j \\
		\sum \limits_{1}^{n} \boldsymbol{O}_j \times \boldsymbol{x}_j \delta \rho_j \\
	\end{array}\right]+O(\delta \boldsymbol{\rho}^2)
\end{aligned}
\end{equation}
The first-order error is the rotation described above, and the change of signed distance from origin of the global coordinate system $\boldsymbol{0}$ to the $xy$-plane of local coordinate system built on panel $n$, as a circuit is taken around the boundary of the hole with folding angles $\boldsymbol{\rho}+\delta \boldsymbol{\rho}$.

For a single vertex, a graphical representation of first-order errors $\varepsilon_1,\varepsilon_2,\varepsilon_3$ is provided in figure \ref{fig: ring}. The consistency constraint is illustrated by a closed torus. When there is a cut, the constraint is released, and the first-order errors are shown by the rotation of cross-section of the torus. Suppose $n=1$, a first-order error $\boldsymbol{\varepsilon}$ would be a rotation $\delta \rho_1$ along direction $\boldsymbol{x}_1$, in (b)(i) and (c)(ii) of figure \ref{fig: ring} we could see that a positive first-order error when $\delta \rho_1>0$ is compatible with a positive folding angle if considering the cut as rotation of panels around an inner crease. 

\subsection{Rigidity matrix for a large rigid origami} \label{subsection: assemble}

Now we consider assembling the derivative for each single-vertex or single-hole in equations \eqref{eq: Jacvertex} and \eqref{eq: Jachole} to a large rigid origami. In the view of programming, information of the crease pattern could be stored in a \textit{incidence matrix} $D$ describing the relationship between inner creases and vertices with a labelling of them. If vertex $i$ is incident to inner crease $j$ and the direction vector goes out from $i$, $D_{ij}=1$; if the direction goes toward $i$, $D_{ij}=-1$, otherwise $D_{ij}=0$. $D$ is a sparse matrix. An example is provided in Section S3.1 of the supplementary material.

\section{Static rigidity} \label{section: statics}

We will now consider the behaviour of a rigid origami when load is applied. First we will introduce a restricted set of external applied loads and internal forces that are work-conjugate to the kinematic quantities mentioned in the previous section, before pointing out how these might be related to more general sets of forces.

The equilibrium analysis starts from the principle of virtual work. A virtual displacement for rigid origami is exactly an arbitrarily small first-order flex $ \boldsymbol{\rho}'$ at a rigidly folded state $(\boldsymbol{\rho}, \lambda')$. We define a \textit{load} $\boldsymbol{l}$ so that the external virtual work done by the load, $\delta W_e$, for any $\boldsymbol{\rho}'$, is given by 
\begin{equation}
\delta W_e = \boldsymbol{l} \boldsymbol{\rho}'.
\end{equation}
Thus the form of the load must be a set of equal and opposite torques applied to the panels on each side of each inner crease, such that positive external virtual work is done by a positive change in folding angle at the crease (remembering that a valley fold corresponds with the positive direction of folding angle). 

\begin{figure}[hbtp]
\noindent \begin{centering}
	\includegraphics[scale=0.6]{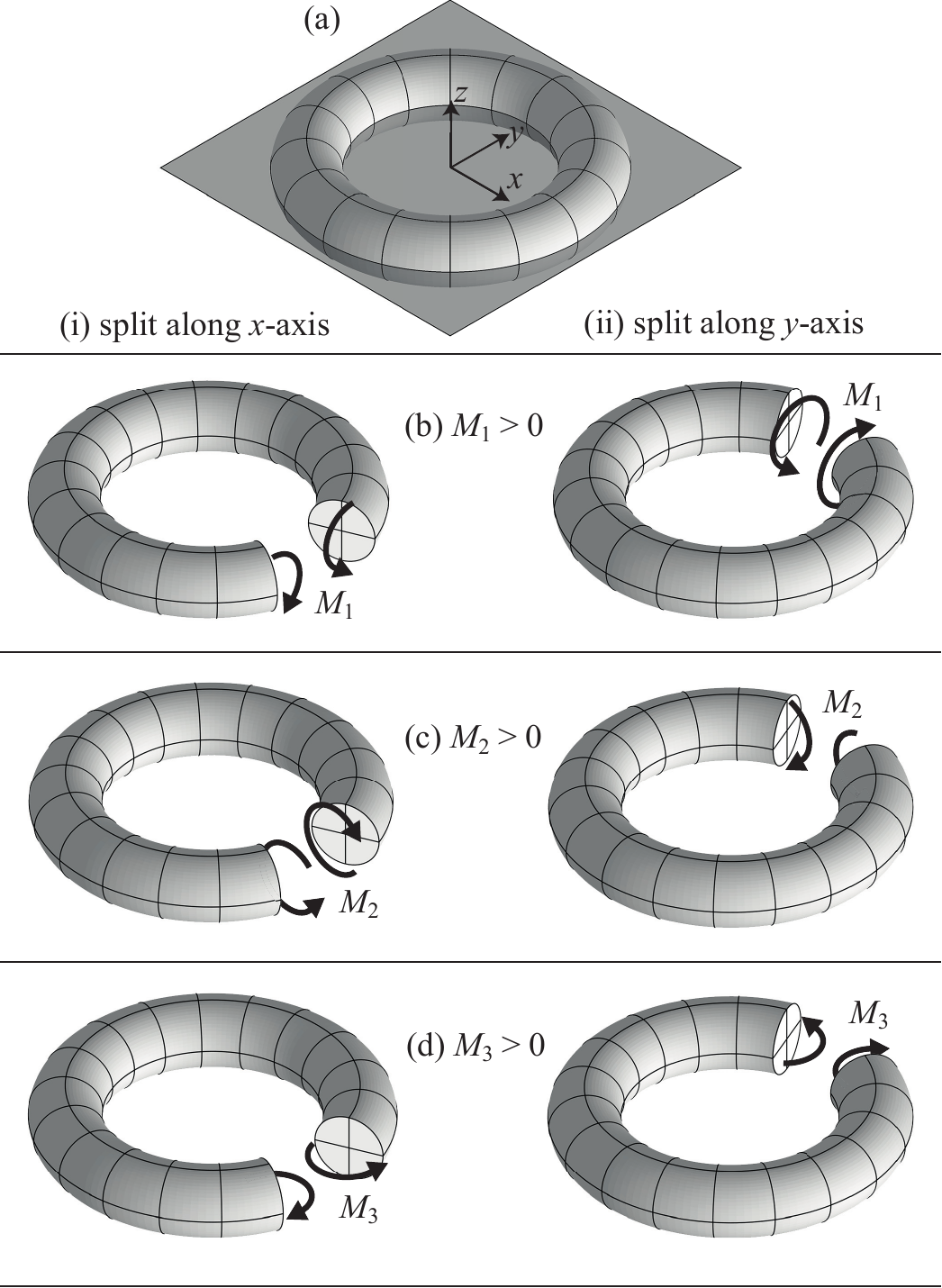}
	\par\end{centering}
\caption{\label{fig: intforces}The internal forces (here, torques) $\{M_1, M_2, M_3\}$ that are work-conjugate to the first-order errors $\{\varepsilon_1, \varepsilon_2, \varepsilon_3\}$ distributed around a torus embedded in a paper, as introduced in figure~\ref{fig: ring}. In each case, the internal force is only shown where there is a cut.  Only in (b)(i) and (c)(ii) can these internal forces be directly applied by a pair of opposite torques at a hinge (along the $x$- and $y$-axes, respectively).  The deformations shown in figure~\ref{fig: ring} are not, \emph{in general}, elastic responses to the forces shown here, although in fact they would correspond to the elastic response of a cut torus that had equal bending and torsional stiffness.}
\end{figure}

Consider also the internal forces that may exist within the rigid panels.  We define the \emph{internal forces} $\boldsymbol{\omega}$ such that the internal work done $\delta W_i$, for any 
first-order error $\boldsymbol{\varepsilon}(\boldsymbol{\rho}')$, is given by 
\begin{equation}
\delta W_i = \boldsymbol{\omega} \boldsymbol{\varepsilon}(\boldsymbol{\rho}') = \boldsymbol{\omega} \dfrac{\mathrm{d} \boldsymbol{A}}{\mathrm{d} \boldsymbol{\rho}} \boldsymbol{\rho}'
\end{equation}
Thus the form of internal force should be an internal torque around each inner vertex; and an internal torque and force around each hole. For a single vertex, a graphical representation of the internal force $\boldsymbol{\omega}$ is provided in figure \ref{fig: intforces}.

The sufficient and necessary condition for equilibrium is $\delta W_i=\delta W_e$, hence we have:
\begin{equation} \label{eq: equilibrium}
\boldsymbol{\omega} \dfrac{\mathrm{d} \boldsymbol{A}}{\mathrm{d} \boldsymbol{\rho}} =\boldsymbol{l}
\end{equation}
For rigid origami, where the error $\boldsymbol{\varepsilon}(\boldsymbol{\rho}')$ is zero,
\begin{equation} 
\delta W_e=\boldsymbol{l}  \boldsymbol{\rho}'= \delta W_i= 0
\end{equation}
For zero load, the \textit{self-stress} $\boldsymbol{\omega}_s$ satisfies:
\begin{equation} \label{eq: self-stress}
\boldsymbol{\omega}_s \dfrac{\mathrm{d} \boldsymbol{A}}{\mathrm{d} \boldsymbol{\rho}}  =\boldsymbol{0}
\end{equation}
Physically, the internal forces could be interpreted as the resistance to deformation 
around each single-vertex or single-hole, and can be revealed only by cutting through the rigid panels.

\begin{rem}
In rigidity theory, $\boldsymbol{\omega}$ is usually referred to as a stress, but we will not use that notation here because of the potential for confusion with the related but different use of the term stress in mechanics. However, we will still use the term self-stress for $\boldsymbol{\omega}_s$, a set of internal forces in equilibrium with zero applied load.
\end{rem} 

Table~\ref{tab: compare} shows the correspondence between a bar-joint framework and the model for rigid origami considered here.

\begin{table} [H]
\begin{center}
	\begin{tabular}{c|cc}
		& Bar-joint framework & Rigid origami  \\
		\hline
		Body & joint & inner crease\\
		Number of freedoms of a body & 3 & 1 \\
		Constraints on freedoms & bars & vertices and holes \\
		Form of the internal force & force/length & torque and force \\
		Form of the external load & forces & equal and opposite torques \\
	\end{tabular}
	\caption{\label{tab: compare}A comparison of the statics of bar-joint framework and rigid origami. In rigid origami we consider that a ``body'' is an inner crease with only 1 freedom, and the constraint is on each ``linkage'' of the rigid origami, specifically, each vertex or hole.}
\end{center}
\end{table}

From equation \eqref{eq: equilibrium} and the structure of rigidity matrix $\mathrm{d} \boldsymbol{A}/\mathrm{d} \boldsymbol{\rho}$, the statement of equilibrium for an inner crease incident to a vertex is that the projection of internal forces on the crease must be equal to the load applied. Consider a single-vertex, an internal force is a torque $\{M_1,M_2,M_3\}$ with 3 components in a global coordinate system. For each inner crease $j~(1 \le j \le N_c)$ with direction vector $\boldsymbol{p}_j$, 
\begin{equation} \label{eq: crease equilibrium 1}
\left[
{\begin{array}{ccc}
		M_1 & M_2 & M_3 
\end{array}} \right] \boldsymbol{p}_j =l_j
\end{equation}
For an inner crease incident to a hole, forces also contribute to the equilibrium.  For a single-hole, the internal force $\{M_1,M_2,M_3,F_1,F_2,F_3\}$ has 6 components of torque and force, which satisfies:
\begin{equation} \label{eq: 11}
\left[
{\begin{array}{ccc}
		M_1 & M_2 & M_3
\end{array}} \right] \boldsymbol{p}_j+ \left[
{\begin{array}{ccc}
		F_1 & F_2 & F_3
\end{array}} \right] (\boldsymbol{O}_j \times \boldsymbol{p}_j) =l_j
\end{equation}
which shows that the torque equilibrium is actually
\begin{equation} \label{eq: crease equilibrium 2}
\left(\left[
{\begin{array}{ccc}
		M_1 & M_2 & M_3
\end{array}} \right] + \left[
{\begin{array}{ccc}
		F_1 & F_2 & F_3
\end{array}} \right] \times \boldsymbol{O}_j \right) \boldsymbol{p}_j =l_j
\end{equation}

Now that we have clarified the form of the internal forces and load, we can consider static rigidity, and its relation to first-order rigidity:
\begin{defn}
A rigid origami $(\boldsymbol{\rho}, \lambda')$ can \textit{resolve} a load $\boldsymbol{l}$ if there is an internal force satisfying equation \eqref{eq: equilibrium}. A rigid origami is \textit{statically rigid} if can resolve every load. A rigid origami is \textit{independent} if there is only zero self-stress. A rigid origami is \textit{isostatic} if first-order rigid and independent.
\end{defn}

\begin{thm}
For a rigid origami $(\boldsymbol{\rho}, \lambda')$ with $N_v$ inner vertices, $N_h$ holes and $N_c$ inner creases, The following statements are equivalent:
\begin{enumerate} [(1)]
	\item $(\boldsymbol{\rho}, \lambda')$ is first-order rigid.
	\item $(\boldsymbol{\rho}, \lambda')$ is statically rigid.
	\item The dimension of the collection of self-stress at $(\boldsymbol{\rho}, \lambda')$ is $3N_v+6N_h-N_c$.
\end{enumerate}
\end{thm}

\begin{proof}
For the rigidity matrix $\mathrm{d} \boldsymbol{A}/\mathrm{d} \boldsymbol{\rho}$, a zero nullspace is equivalent to a full image. The rank of its left nullspace when the nullspace is zero is $3N_v+6N_h-N_c$.
\end{proof}

Some examples showing the calculation of internal forces and states of self-stress are given in Section S3.2 of the supplementary material.

\begin{rem}
If any set of self-equilibrating forces and torques are applied (either discrete or continuous, which we describe as a \textit{general load}), there would be a unique decomposition to the form of load defined above, i.e. opposite torques applied on adjacent panels around each inner crease. The rest of a general load could always be carried by a rigid origami. To examine this, a possible way is to replace the rigid origami by a corresponding double-coning bar-joint framework. Such framework is generated by replacing the boundary of each panel by a series of bars and joints. Then adding two out-of-plane vertices on different side of the panel, and joining the two vertices to each of the vertex on the panel with a bar.
\end{rem}

\section{Prestress stability}  \label{section: prestress stability}

In this section we consider rigid origami that are not first-order rigid, but are rigid, and elucidate how the stability of these structures is changed when prestress added.  To do that, we will describe an energy function $U$ that gives the potential energy stored in the paper. In fact, for our purposes the energy function can be fairly general in its form, but it can also be given in a quite physical way.  We will see that the first differential of $U$ with respect to the folding angles naturally gives a state of self-stress for the paper, and the second differential naturally leads to the stiffness, and hence stability.

\subsection{Energy, stiffness and stability}

\begin{defn}
The potential energy $U$ stored in a rigid origami only depends on the error of independent consistency constraint around the $N_v$ inner vertices and $N_h$ holes, $\boldsymbol{A} \in \mathbb{R}^{3N_v+6N_h}$, and satisfies
\begin{equation} 
	U(\boldsymbol{0})=0, \quad U(\boldsymbol{A})>0 ~~ \mathrm{if} ~~ \boldsymbol{A} \ne \boldsymbol{0}
\end{equation} 
We require $U$ to have continuous second-order derivative, so that we can define the matrix $B$ as
\begin{equation} \label{eq: energy}
	B 		= \dfrac{\mathrm{d}^2 U} 
	{\mathrm{d} \boldsymbol{A}^2}, 
	\quad B_{il} 	= \dfrac{\partial^2 U} 
	{\partial A_i \, \partial A_l}.
\end{equation}
The size of $B$ is $(3N_v+6N_h) \times (3N_v+6N_h)$, which is assumed to be positive definite. (Note that we are using $i$ and $l$ as subscripts corresponding to error components; later $j$ and $k$ will be used as subscripts for folding angles.)
\end{defn}

\begin{rem}
We could consider a less general energy function than that given in \eqref{eq: energy}, where the energy is the sum of the energy stored around each inner vertex or hole, in which case the matrix $\boldsymbol{B}$ will be block-diagonal, with one block per inner vertex or hole.  Or we might wish to consider that the energy stored by each misfit error $A_i$ is independent, so that $\boldsymbol{B}$ is diagonal.  For a particularly simple choice, we could define $E = \sum \frac{1}{2} g_i A_i^2$, so that $\boldsymbol{B}$ would be diagonal and constant, with $B_{ii} = g_i$.  All of these choices might impact the physical behaviour of the system under load, but will not affect the definition of prestress stiffness below.
\end{rem}

Next we will consider the equilibrium of a rigid origami from an energy viewpoint, and judge whether it is stable. In general, suppose the rigid origami is in a conservative force field with potential $V(\boldsymbol{\rho})$, then the total energy could be written as 
\begin{equation}
E=U+V
\end{equation}
The partial derivative of $E$ with respect to a folding angle $\rho_j$ is ($1 \le i \le 3N_v+6N_h$, $1 \le j \le N_c$)
\begin{equation} 
\begin{gathered}
	\dfrac{\partial E}{\partial \rho_j}= \dfrac{\partial U}{\partial A_i} \dfrac{\partial A_i}{\partial \rho_j}+\dfrac{\partial V}{\partial \rho_j} \\
\end{gathered}
\end{equation}
which can be written in a more compact form,
\begin{equation} \label{eq: first-order energy}
\dfrac{\mathrm{d} E}{\mathrm{d} \boldsymbol{\rho}}= \dfrac{\mathrm{d} U} {\mathrm{d} \boldsymbol{A}}  \dfrac{\mathrm{d} \boldsymbol{A}}{\mathrm{d} \boldsymbol{\rho}}+\dfrac{\mathrm{d} V}{\mathrm{d} \boldsymbol{\rho}}
\end{equation}
The equilibrium condition is then
\begin{equation} 
\dfrac{\mathrm{d} U} {\mathrm{d} \boldsymbol{A}} \dfrac{\mathrm{d} \boldsymbol{A}}{\mathrm{d} \boldsymbol{\rho}}+\dfrac{\mathrm{d} V}{\mathrm{d} \boldsymbol{\rho}}= \boldsymbol{0}
\end{equation}
Since 
\begin{equation} 
\boldsymbol{l}=-\dfrac{\mathrm{d} V}{\mathrm{d} \boldsymbol{\rho}}
\end{equation}
the above condition is exactly equation \eqref{eq: equilibrium}, which shows that the first-order derivative of the energy function $U$ is an internal force $\boldsymbol{\omega}$. When ${\mathrm{d} V}/{\mathrm{d} \boldsymbol{\rho}}=\boldsymbol{0}$, the first-order derivative of $U$ is a self-stress.

To consider stability of the equilibrium we have to consider the second differential, the Hessian of energy, ($1 \le i,l \le 3N_v+6N_h, 1 \le j,k \le N_c$)
\begin{equation} 
\begin{gathered}
	\dfrac{\partial^2 E}{\partial \rho_k \partial \rho_j}=\dfrac{\partial^2 U}{\partial A_l \partial A_i} \dfrac{\partial A_l}{\partial \rho_k} \dfrac{\partial A_i}{\partial \rho_j}+\dfrac{\partial U}{\partial A_i} \dfrac{\partial^2 A_i}{\partial \rho_k \partial \rho_j}+\dfrac{\partial^2 V}{\partial \rho_k \partial \rho_j}\\
\end{gathered}
\end{equation}
which can be written in a compact form,
\begin{equation} \label{eq: second-order energy}
\dfrac{\mathrm{d}^2 E}{\mathrm{d} \boldsymbol{\rho}^2}= \dfrac{\mathrm{d} \boldsymbol{A}}{\mathrm{d} \boldsymbol{\rho}}^T \dfrac{{\mathrm{d}^2 U}}{\mathrm{d} \boldsymbol{A}^2} \dfrac{\mathrm{d} \boldsymbol{A}}{\mathrm{d} \boldsymbol{\rho}}+\dfrac{\mathrm{d} U} {\mathrm{d} \boldsymbol{A}} \dfrac{{\mathrm{d}^2 \boldsymbol{A}}}{\mathrm{d} \boldsymbol{\rho}^2}-\dfrac{\mathrm{d} \boldsymbol{l}}{\mathrm{d} \boldsymbol{\rho}}
\end{equation}
The condition for stability is that the total energy at a rigidly folded state reaches a strict local minimum, and a sufficient condition is the second-order differential of the total energy is positive definite. The second-order derivative ${\mathrm{d}^2 \boldsymbol{A}}/\mathrm{d} \boldsymbol{\rho}^2$ is also called the Hessian of the independent consistency constraint $\boldsymbol{A}(\boldsymbol{\rho})$, an order 3 tensor with size $(3N_v+6N_h) \times N_c \times N_c$, which could be written in an explicit form as provided in next subsection.

If there is a perturbation of folding angle $\delta \boldsymbol{\rho}$ around a rigidly folded state $(\boldsymbol{\rho},\lambda')$, the increase of total potential energy in second-order will be 
\begin{equation}
\delta E=\dfrac{1}{2} \delta \boldsymbol{\rho}^T \dfrac{\mathrm{d}^2 E}{\mathrm{d} \boldsymbol{\rho}^2} \delta \boldsymbol{\rho}
\end{equation}
and the restoring force will be
\begin{equation}
F=-\dfrac{\partial \delta E}{\partial \delta \boldsymbol{\rho}}=-\dfrac{\mathrm{d}^2 E}{\mathrm{d} \boldsymbol{\rho}^2} \delta \boldsymbol{\rho}
\end{equation}
The above derivation shows how $\mathrm{d}^2 E/\mathrm{d} \boldsymbol{\rho}^2$ works as the stiffness of the rigid origami system. However, if $\delta E=0$ for a perturbation $\delta \boldsymbol{\rho}$, for this direction we might need higher order information of energy to determine the stability.

In this section we will discuss the prestress stability first, assuming there is no load (${\mathrm{d} V}/{\mathrm{d} \boldsymbol{\rho}}=\boldsymbol{0}$).
\begin{defn} \label{defn: prestress stability}
A rigid origami $(\boldsymbol{\rho},\lambda')$ with $N_v$ inner vertices, $N_c$ inner creases and $N_h$ holes is \textit{prestress stable} if there is a positive-definite matrix $B$ with size ($3N_v+6N_h) \times (3N_v+6N_h$), and a vector $\boldsymbol{\omega}_s \in \mathbb{R}^{3N_v+6N_h}$ such that
\begin{equation}
	\boldsymbol{\omega}_s \dfrac{\mathrm{d} \boldsymbol{A}}{\mathrm{d} \boldsymbol{\rho}}=\boldsymbol{0}
\end{equation}
and
\begin{equation} \label{eq: tsm}
	K=\dfrac{\mathrm{d} \boldsymbol{A}}{\mathrm{d} \boldsymbol{\rho}}^T  B \dfrac{\mathrm{d} \boldsymbol{A}}{\mathrm{d} \boldsymbol{\rho}}+\boldsymbol{\omega}_s \dfrac{{\mathrm{d}^2 \boldsymbol{A}}}{\mathrm{d} \boldsymbol{\rho}^2}
\end{equation}
is positive-definite.  
\end{defn}
Physically, $B$ is the \textit{local elasticity matrix}, which is the Hessian of the predefined energy function. $K$ is the \textit{tangent stiffness matrix} or \textit{total stiffness matrix}. $\boldsymbol{\omega}_s \cdot \mathrm{d}^2 \boldsymbol{A}/\mathrm{d} \boldsymbol{\rho}^2$ is called the \textit{stress matrix}. We say this self-stress $\boldsymbol{\omega}_s$ \textit{stabilizes} a rigid origami if it leads to a positive-definite stiffness $K$. 

\subsection{Hessian for a single-vertex or hole}

We will show how to derive the Hessian ${\mathrm{d}^2 \boldsymbol{A}}/\mathrm{d} \boldsymbol{\rho}^2$ for a large rigid origami followed by writing ${\mathrm{d}^2 \boldsymbol{A}}/\mathrm{d} \boldsymbol{\rho}^2$ for its restriction on a degree-$n$ single-vertex or hole. Consider the second-order derivative of equations \eqref{eq: independent consistency 1} and \eqref{eq: independent consistency 2}. 

When $1 \le k \le j \le n$,
\begin{equation} \label{eq: Hessian 1}
\dfrac{\partial^2 }{\partial \rho_k \partial \rho_j} 	\left[\begin{array}{c}
	A_1 \\
	A_2 \\
	A_3
\end{array}\right]= \left[\begin{array}{c}
	x_{2k}x_{3j} \\
	x_{3k}x_{1j} \\
	x_{1k}x_{2j}
\end{array}\right]
\end{equation}
where $\boldsymbol{x}_j=[x_{1j};x_{2j};x_{3j}]$ is the direction vector of the inner crease $\rho_j$ measured in a global coordinate system, pointing away from this vertex.
\begin{equation} \label{eq: Hessian 2}
\dfrac{\partial^2 }{\partial \rho_k \partial \rho_j} \left[\begin{array}{c}
	A_1 \\
	A_2 \\
	A_3 \\
	A_4 \\
	A_5 \\
	A_6
\end{array}\right] = \left[
\begin{array}{c}
	x_{2k}x_{3j} \\
	x_{3k}x_{1j} \\
	x_{1k}x_{2j} \\
	\\
	\boldsymbol{x}_k \times (\boldsymbol{O}_j \times \boldsymbol{x}_j) \\
	\\ 
\end{array}
\right]
\end{equation}
where $\boldsymbol{x}_j=[x_{1j};x_{2j};x_{3j}]$ is the direction vector of the inner crease $\rho_j$ measured in a global coordinate system, pointing away from this hole. $\boldsymbol{O}_j$ is the position of vertex on the hole incident to $\rho_j$ measured in the global coordinate system (figure \ref{fig: single vertex and hole}). When $1 \le j < k \le n$, swap $k$ and $j$. The derivation is provided in Section S2 of the supplementary material.

\subsection{Hessian for a large rigid origami}

Now we consider assembling the second-order derivative for each single-vertex or single-hole in equations \eqref{eq: Hessian 1} and \eqref{eq: Hessian 2} to a large rigid origami. The Hessian ${\mathrm{d}^2 \boldsymbol{A}}/\mathrm{d} \boldsymbol{\rho}^2$ is an order 3 tensor with size $(3N_v+6N_h) \times N_c \times N_c$, which is the collection of a series of sparse matrices ${\mathrm{d}^2 A_i}/\mathrm{d} \boldsymbol{\rho}^2$. The size for each of them is $N_c \times N_c$.

Consider the incidence matrix for vertices $D_{\mathrm{vertex}}$, row $i(1 \le i \le N_v)$ corresponds with three components $[A_{3i-2}; A_{3i-1}; A_{3i}]$. For folding angle $\rho_j$ where the direction vector is $\boldsymbol{p}_j$ in a global coordinate system, if $D_{ij}=1$, $\boldsymbol{x}_j=\boldsymbol{p}_j$, otherwise if $D_{ij}=-1$, $\boldsymbol{x}_j=-\boldsymbol{p}_j$, then apply equation \eqref{eq: Hessian 1} we could obtain $3N_v$ matrices ${\mathrm{d}^2 A_i}/\mathrm{d} \boldsymbol{\rho}^2$ for every single-vertex.

Next consider the incidence matrix for holes $D_{\mathrm{hole}}$, row $i(1 \le i \le N_h)$ corresponds with six components $[A_{3N_v+6i-5}; A_{3N_v+6i-4}; A_{3N_v+6i-3}; A_{3N_v+6i-2}; A_{3N_v+6i-1}; A_{3N_v+6i}]$. For folding angle $\rho_j$, if $D_{ij}=1$, $\boldsymbol{x}_j=\boldsymbol{p}_j$, otherwise if $D_{ij}=-1$, $\boldsymbol{x}_j=-\boldsymbol{p}_j$, then apply equation \eqref{eq: Hessian 2} we could obtain $6N_h$ matrices ${\mathrm{d}^2 A_i}/\mathrm{d} \boldsymbol{\rho}^2$ for every single-hole.

\subsection{Reducing the calculation}

From definition \ref{defn: prestress stability}, if a rigid origami is first-order rigid-foldable and there is no prestress applied, it is not stable. Hence for a given configuration, an important question is to find the collection of self-stress to stabilize a rigid origami. The next proposition provides a simpler way to judge the prestress stability.

\begin{prop} \label{prop: 5}
Concerning the positive (semi)-definiteness of the stiffness matrix for a first-order rigid-foldable rigid origami.
\begin{enumerate} [(1)]
	\item The matrix $[(\mathrm{d} \boldsymbol{A}/\mathrm{d} \boldsymbol{\rho})^T \cdot B \cdot \mathrm{d} \boldsymbol{A}/\mathrm{d} \boldsymbol{\rho}]$ is positive semi-definite.  Its quadratic and linear nullspace are both the nullspace of $\mathrm{d} \boldsymbol{A}/\mathrm{d} \boldsymbol{\rho}$.
	\item A rigid origami $(\boldsymbol{\rho},\lambda')$ is prestress stable if and only if there exists a self-stress $\boldsymbol{\omega}_s \in \mathbb{R}^{3N_v+6N_h}$ such that $\boldsymbol{\omega}_s \cdot \mathrm{d}^2 \boldsymbol{A}/\mathrm{d} \boldsymbol{\rho}^2$ is positive-definite when restricted to the nullspace of $\mathrm{d} \boldsymbol{A}/\mathrm{d} \boldsymbol{\rho}$. 
	\item A rigid origami $(\boldsymbol{\rho},\lambda')$ is prestress stable if and only if there exists a self-stress $\boldsymbol{\omega}_s \in \mathbb{R}^{3N_v+6N_h}$ such that all the eigenvalues of $[\rho']^T \cdot \boldsymbol{\omega}_s \cdot \mathrm{d}^2 \boldsymbol{A}/\mathrm{d} \boldsymbol{\rho}^2 \cdot [\rho']$ are positive, where $[\rho']$ is the collection of a basis of first-order flex.
	\begin{equation} \label{eq: first-order basis}
		[\rho']= \left[\begin{array}{cccc}
			\boldsymbol{\rho}_1' & \boldsymbol{\rho}_2' & \cdots & \boldsymbol{\rho}_{N_c-\mathrm{rank}(\mathrm{d} \boldsymbol{A}/\mathrm{d} \boldsymbol{\rho})}'
		\end{array} \right]
	\end{equation}
\end{enumerate}
\end{prop}

\begin{proof} 
Statement (1): For any perturbation of folding angles $\delta \boldsymbol{\rho}$, consider the quadratic form, if $\delta \boldsymbol{\rho}$ is not a first-order flex, since $B$ is assumed to be positive-definite, $\delta \boldsymbol{\rho}^T \cdot  (\mathrm{d} \boldsymbol{A}/\mathrm{d} \boldsymbol{\rho})^T \cdot B \cdot \mathrm{d} \boldsymbol{A}/\mathrm{d} \boldsymbol{\rho} \cdot \delta \boldsymbol{\rho}>0$, therefore $[(\mathrm{d} \boldsymbol{A}/\mathrm{d} \boldsymbol{\rho})^T \cdot B \cdot \mathrm{d} \boldsymbol{A}/\mathrm{d} \boldsymbol{\rho}]$ is positive semi-definite, and the quadratic nullspace of $[(\mathrm{d} \boldsymbol{A}/\mathrm{d} \boldsymbol{\rho})^T \cdot B \cdot \mathrm{d} \boldsymbol{A}/\mathrm{d} \boldsymbol{\rho}]$ is contained in the collection of first-order flexes given by the nullspace of $\mathrm{d} \boldsymbol{A}/\mathrm{d} \boldsymbol{\rho}$. Since the nullspace of $\mathrm{d} \boldsymbol{A}/\mathrm{d} \boldsymbol{\rho}$ is also contained in the quadratic nullspace of $[(\mathrm{d} \boldsymbol{A}/\mathrm{d} \boldsymbol{\rho})^T \cdot B \cdot \mathrm{d} \boldsymbol{A}/\mathrm{d} \boldsymbol{\rho}]$, the statement holds. Then consider the linear form. For any $\delta \boldsymbol{\rho}$, if $(\mathrm{d} \boldsymbol{A}/\mathrm{d} \boldsymbol{\rho})^T \cdot B \cdot \mathrm{d} \boldsymbol{A}/\mathrm{d} \boldsymbol{\rho} \cdot \delta \boldsymbol{\rho}=\boldsymbol{0}$, $\delta \boldsymbol{\rho}^T \cdot (\mathrm{d} \boldsymbol{A}/\mathrm{d} \boldsymbol{\rho})^T \cdot B \cdot \mathrm{d} \boldsymbol{A}/\mathrm{d} \boldsymbol{\rho} \cdot \delta \boldsymbol{\rho}=0$, therefore the linear nullspace is contained in the quadratic nullspace. If $\delta \boldsymbol{\rho}$ is in the quadratic nullspace, it is a first-order flex, therefore being an element of the linear nullspace, hence the statement holds.

Statement (2): Necessity: if a rigid origami is prestress stable, the quadratic form of a first-order flex should be greater than $0$, hence $\boldsymbol{\omega}_s \cdot {\mathrm{d}^2 \boldsymbol{A}}/\mathrm{d} \boldsymbol{\rho}^2$ is positive-definite when restricted to the nullspace of  $\mathrm{d} \boldsymbol{A}/\mathrm{d} \boldsymbol{\rho}$.

Sufficiency: We will show that, if there exists a self-stress $\boldsymbol{\omega}_s$ such that $\boldsymbol{\omega}_s \cdot {\mathrm{d}^2 \boldsymbol{A}}/\mathrm{d} \boldsymbol{\rho}^2$ is positive definite when restricted to the nullspace of $\mathrm{d} \boldsymbol{A}/\mathrm{d} \boldsymbol{\rho}$, by choosing a sufficiently large $k$, $K=(\mathrm{d} \boldsymbol{A}/\mathrm{d} \boldsymbol{\rho})^T \cdot B \cdot \mathrm{d} \boldsymbol{A}/\mathrm{d} \boldsymbol{\rho}+k\boldsymbol{\omega}_s \cdot {\mathrm{d}^2 \boldsymbol{A}}/\mathrm{d} \boldsymbol{\rho}^2$ would be positive-definite. 

For any perturbation of folding angles $\delta \boldsymbol{\rho}$, if $\delta \boldsymbol{\rho}$ is a first-order flex, for any $k>0$, $\delta \boldsymbol{\rho}^T K \delta \boldsymbol{\rho} > 0$. If $\delta \boldsymbol{\rho}$ is not a first-order flex, suppose $||\delta \boldsymbol{\rho}||=1$.  Since this set is compact, there exists $\varepsilon>0$ s.t. $\delta \boldsymbol{\rho}^T \cdot (\mathrm{d} \boldsymbol{A}/\mathrm{d} \boldsymbol{\rho})^T \cdot B \cdot \mathrm{d} \boldsymbol{A}/\mathrm{d} \boldsymbol{\rho} \cdot \boldsymbol{\rho} \ge \varepsilon$ and we know $\delta \boldsymbol{\rho}^T \cdot \boldsymbol{\omega}_s \cdot {\mathrm{d}^2 \boldsymbol{A}}/\mathrm{d} \boldsymbol{\rho}^2 \cdot \delta \boldsymbol{\rho} \ge -||\boldsymbol{\omega}_s \cdot {\mathrm{d}^2 \boldsymbol{A}}/\mathrm{d} \boldsymbol{\rho}^2||$, then we can choose $0<k<\varepsilon/||\boldsymbol{\omega}_s \cdot {\mathrm{d}^2 \boldsymbol{A}}/\mathrm{d} \boldsymbol{\rho}^2||$ s.t. $\delta \boldsymbol{\rho}^T K \delta \boldsymbol{\rho} > 0$. Further, consider $||\delta \boldsymbol{\rho}|| \neq 1$, choosing the same $k$, $\delta \boldsymbol{\rho}^T K \delta \boldsymbol{\rho}/||\delta \boldsymbol{\rho}||^2 > 0$.  

Statement (3): Since a first-order flex $\boldsymbol{\rho}'$ could be written as
and the collection of coefficients as
\begin{equation}
	\boldsymbol{\rho}'=[\rho']\boldsymbol{a} = \left[\begin{array}{cccc}
		\boldsymbol{\rho}_1' & \boldsymbol{\rho}_2' & \cdots & \boldsymbol{\rho}_{N_c-\mathrm{rank}(\mathrm{d} \boldsymbol{A}/\mathrm{d} \boldsymbol{\rho})}'
	\end{array} \right] \left[\begin{array}{c}
		a_1\\ 
		a_2 \\ 
		\vdots \\
		a_{N_c-\mathrm{rank}(\mathrm{d} \boldsymbol{A}/\mathrm{d} \boldsymbol{\rho})}
	\end{array} \right]
\end{equation}
this statement will hold from (2).
\end{proof} 

Calculation of some specific examples that are prestress stable but not first-order rigid are provided in S3.3 of the supplementary material. There is a general result for a rigid and planar single-vertex. The proof is given in S4 of the supplementary material.

\begin{prop} \label{prop: single-vertex prestress stability}
A rigid planar single-vertex is prestress stable.
\end{prop}

Further, any triangulated convex polyhedral surface (possibly with some properly placed holes) is prestress stable \cite{connelly_prestress_2017}. This is drawn from the analysis of a ``spider tensegrity'', and would be applicable to rigid origami. Note that from Dehn's result \cite{dehn_uber_1916} on Cauchy's rigidity theory \cite{cauchy_sur_1813}, a strictly convex polyhedral surface is first-order rigid.

\subsection{Considering external load}
When there is a load $\boldsymbol{l}(\boldsymbol{\rho})$ applied on a rigid origami, the above theory on stability could be modified as below.

\begin{defn} \label{defn: stability}
A rigid origami $(\boldsymbol{\rho},\lambda')$ with $N_v$ inner vertices, $N_c$ inner creases and $N_h$ holes is stable under load $\boldsymbol{l}(\boldsymbol{\rho})$ if there is a positive-definite matrix $B$ with size ($3N_v+6N_h) \times (3N_v+6N_h$), and a vector $\boldsymbol{\omega} \in \mathbb{R}^{3N_v+6N_h}$ such that
\begin{equation}
	\boldsymbol{\omega} \dfrac{\mathrm{d} \boldsymbol{A}}{\mathrm{d} \boldsymbol{\rho}}=\boldsymbol{l}
\end{equation}
and
\begin{equation} \label{eq: stiffness}
	K=\dfrac{\mathrm{d} \boldsymbol{A}}{\mathrm{d} \boldsymbol{\rho}}^T  B \dfrac{\mathrm{d} \boldsymbol{A}}{\mathrm{d} \boldsymbol{\rho}}+\boldsymbol{\omega} \dfrac{\mathrm{d}^2 \boldsymbol{A}}{\mathrm{d} \boldsymbol{\rho}^2}-\dfrac{\mathrm{d} \boldsymbol{l}}{\mathrm{d} \boldsymbol{\rho}}
\end{equation}
is positive-definite.
\end{defn}

\begin{prop} 
A rigid origami $(\boldsymbol{\rho}, \lambda')$ is stable under load $\boldsymbol{l}(\boldsymbol{\rho})$ if and only if there exists a stress $\boldsymbol{\omega} \in \mathbb{R}^{3N_v+6N_h}$ such that $\boldsymbol{\omega} \cdot \mathrm{d}^2 \boldsymbol{A}/\mathrm{d} \boldsymbol{\rho}^2$ is positive-definite when restricted to the nullspace of  $\mathrm{d} \boldsymbol{A}/\mathrm{d} \boldsymbol{\rho}$. Equivalently, if and only if there exists a stress $\boldsymbol{\omega} \in \mathbb{R}^{3N_v+6N_h}$ such that all the eigenvalues of $[\rho']^T \cdot \boldsymbol{\omega} \cdot \mathrm{d}^2 \boldsymbol{A}/\mathrm{d} \boldsymbol{\rho}^2 \cdot [\rho']$ are positive, where $[\rho']$ is the collection of a basis of first-order flex in equation \eqref{eq: first-order basis}.
\end{prop}

\begin{proof}
From statement (2) of proposition \ref{prop: 5}, $(\boldsymbol{\rho}, \lambda')$ is stable if and only if there exists a stress $\boldsymbol{\omega} \in \mathbb{R}^{3N_v+6N_h}$ such that $\boldsymbol{\omega} \cdot \mathrm{d}^2 \boldsymbol{A}/\mathrm{d} \boldsymbol{\rho}^2-\mathrm{d} \boldsymbol{l}/\mathrm{d} \boldsymbol{\rho}$ is positive-definite when restricted to the nullspace of  $\mathrm{d} \boldsymbol{A}/\mathrm{d} \boldsymbol{\rho}$. Since a first-order flex $\boldsymbol{\rho}'$ is orthogonal to $\boldsymbol{l}$, the quadratic form $\boldsymbol{\rho}' \cdot \mathrm{d} \boldsymbol{l}/\mathrm{d} \boldsymbol{\rho} \cdot \boldsymbol{\rho}=0$.
\end{proof}

\begin{rem}
In equation \eqref{eq: stiffness}, the first term $(\mathrm{d} \boldsymbol{A}/\mathrm{d} \boldsymbol{\rho})^T \cdot B \cdot \mathrm{d} \boldsymbol{A}/\mathrm{d} \boldsymbol{\rho}$ could be interpreted as the material part of the stiffness matrix, which only relates to how the potential energy stored in a rigid origami is defined, and assumed to be semi positive-definite with nullspace as the collection of first-order flex. The second term $\boldsymbol{\omega} \cdot {\mathrm{d}^2 \boldsymbol{A}}/\mathrm{d} \boldsymbol{\rho}^2$ shows how a load could possibly enhance or reduce the stiffness. Further, the third term $-\mathrm{d} \boldsymbol{l}/\mathrm{d} \boldsymbol{\rho}$ will change the restoring force but has no effect on the stability.
\end{rem}

An example showing how load would affect the stability is also given in Section S3.3 of the supplementary material.

\section{Second-order Rigidity} \label{section: second order rigidity}

In this section we will discuss the second-order rigidity and show its link with prestress stability. For a rigid origami $(\boldsymbol{\rho},\lambda')$, a first-order flex $\boldsymbol{\rho}'$ is obtained by differentiating the independent consistency constraint $\boldsymbol{A}(\boldsymbol{\rho})=\boldsymbol{0}$. Similarly, a second-order flex $(\boldsymbol{\rho}',\boldsymbol{\rho}'')$ satisfies the condition from differentiating the consistency constraint twice.

\begin{defn} \label{defn: second order}
For a rigid origami $(\boldsymbol{\rho},\lambda')$ with $N_v$ inner vertices, $N_c$ inner creases and $N_h$ holes, a \textit{second-order flex} $(\boldsymbol{\rho}',\boldsymbol{\rho}'') \in (\mathbb{R}^{N_c},\mathbb{R}^{N_c})$ satisfies ($1 \le i \le 3N_v+6N_h, 1 \le j,k \le N_c$)
\begin{equation} \label{eq: second-order}
	\begin{gathered}
		\dfrac{\partial A_i}{\partial \rho_j}\rho_j'=0 \\
		\dfrac{\partial^2 A_i}{\partial \rho_j \partial \rho_k}\rho_j'\rho_k'+	\dfrac{\partial A_i}{\partial \rho_j}\rho_j''=0
	\end{gathered}
\end{equation}
If written in a compact form,
\begin{equation} 
	\begin{gathered}
		\dfrac{\mathrm{d} \boldsymbol{A}}{\mathrm{d} \boldsymbol{\rho}} \boldsymbol{\rho}'=\boldsymbol{0} \\
		\boldsymbol{\rho}'^T  \dfrac{\mathrm{d}^2 \boldsymbol{A}}{\mathrm{d} \boldsymbol{\rho}^2}  \boldsymbol{\rho}' + \dfrac{\mathrm{d} \boldsymbol{A}}{\mathrm{d} \boldsymbol{\rho}}  \boldsymbol{\rho}''=\boldsymbol{0}
	\end{gathered}
\end{equation}
A second-order flex with $\boldsymbol{\rho}'=\boldsymbol{0}$ is called \textit{trivial}, otherwise \textit{non-trivial}. If there is only trivial second-order flex, this rigid origami is \textit{second-order rigid}, otherwise \textit{second-order rigid-foldable}.
\end{defn}

\begin{prop} \label{prop: 6}
Some statements concerning the second-order rigidity and prestress stability.
\begin{enumerate}[(1)]
	\item If $(\boldsymbol{\rho}', \boldsymbol{\rho}'')$ is a second-order flex and $\boldsymbol{\rho}_0'$ is a first-order flex, $(\boldsymbol{\rho}', \boldsymbol{\rho}''+\boldsymbol{\rho}_0')$ is also a second-order flex.
	\item A first-order flex $\boldsymbol{\rho}'$ can be extended to a second-order flex $\boldsymbol{\rho}''$ if and only if for all self-stress $\boldsymbol{\omega}_s$,  $\boldsymbol{\rho}'^T  [\boldsymbol{\omega}_s \cdot \mathrm{d}^2 \boldsymbol{A}/\mathrm{d} \boldsymbol{\rho}^2]\, \boldsymbol{\rho}'= 0$.
	\item A rigid origami is second-order rigid if and only if for any first-order flex $\boldsymbol{\rho}'$ there is a self-stress $\boldsymbol{\omega}_s(\boldsymbol{\rho}')$ s.t.,  $\boldsymbol{\rho}'^T [\boldsymbol{\omega}_s \cdot \mathrm{d}^2 \boldsymbol{A}/\mathrm{d} \boldsymbol{\rho}^2]\, \boldsymbol{\rho}' > 0$.
	\item A rigid origami is second-order rigid if and only if the intersection of quadratic nullspace of all $[\rho']^T \cdot \boldsymbol{\omega}_i \cdot \mathrm{d}^2 \boldsymbol{A}/\mathrm{d} \boldsymbol{\rho}^2 \cdot [\rho']$ is $\boldsymbol{0}$. Here $\{\boldsymbol{\omega}_i\}$ is a base of self-stress $(1 \le i \le 3N_v+6N_h-\mathrm{rank}(\mathrm{d} \boldsymbol{A}/\mathrm{d} \boldsymbol{\rho}))$.
	\item A rigid origami is prestress stable if second-order rigid and $\mathrm{rank}(\mathrm{d} \boldsymbol{A}/\mathrm{d} \boldsymbol{\rho})=N_c-1$ or $\mathrm{rank}(\mathrm{d} \boldsymbol{A}/\mathrm{d} \boldsymbol{\rho})=3N_v+6N_h-1$.
\end{enumerate}
\end{prop}

\begin{proof}~

Statement (1) could be verified directly from definition \ref{defn: second order}.

Statement (2): a first-order flex can be extended to a second-order flex if and only if there exists a solution for the linear system below:
\begin{equation}
	\begin{gathered}
		\dfrac{\mathrm{d} \boldsymbol{A}}{\mathrm{d} \boldsymbol{\rho}} \boldsymbol{\rho}''= -\boldsymbol{\rho}'^T \dfrac{\mathrm{d}^2 \boldsymbol{A}}{\mathrm{d} \boldsymbol{\rho}^2} \boldsymbol{\rho}'
	\end{gathered}
\end{equation}
i.e.\ that the vector $(\boldsymbol{\rho}'^T \cdot \mathrm{d}^2 \boldsymbol{A}/\mathrm{d} \boldsymbol{\rho}^2 \cdot \boldsymbol{\rho}')$ lies in the column space of the matrix $\mathrm{d} \boldsymbol{A}/\mathrm{d} \boldsymbol{\rho}$.  Any self-stress $\boldsymbol{\omega}_s$ lies in the left nullspace (the orthogonal complement of the column space) of $\mathrm{d} \boldsymbol{A}/\mathrm{d} \boldsymbol{\rho}$, and hence $\boldsymbol{\omega}_s (\boldsymbol{\rho}'^T   \cdot \mathrm{d}^2 \boldsymbol{A}/\mathrm{d} \boldsymbol{\rho}^2 \cdot \boldsymbol{\rho}')= 0$.  The order of the first two terms in the expression can be swapped without affecting the outcome, and hence the statement is proved.

Statement (3): we know from the inverse negative of statement (2) that a rigid origami is second-order rigid if and only if, for any first order-flex $\boldsymbol{\rho}'$, there is a self-stress $\boldsymbol{\omega}_s(\boldsymbol{\rho}')$ such that  $\boldsymbol{\rho}'^T [\boldsymbol{\omega}_s \cdot \mathrm{d}^2 \boldsymbol{A}/\mathrm{d} \boldsymbol{\rho}^2]\, \boldsymbol{\rho}' \neq 0$.  Either this quadratic form is positive, or can be made positive by replacing $\boldsymbol{\omega}_s$ with $-\boldsymbol{\omega}_s$.

Statement (4): If a first-order flex $\boldsymbol{\rho}'=[\rho'] \boldsymbol{a}$ can be extended to a second-order flex, $\boldsymbol{a}$ should be a quadratic root for every $[\rho']^T \cdot \boldsymbol{\omega}_i \cdot \mathrm{d}^2 \boldsymbol{A}/\mathrm{d} \boldsymbol{\rho}^2 \cdot [\rho']$ such that $\boldsymbol{a}^T [\rho']^T \cdot \boldsymbol{\omega}_i \cdot \mathrm{d}^2 \boldsymbol{A}/\mathrm{d} \boldsymbol{\rho}^2 \cdot [\rho'] \boldsymbol{a}=0$, which leads to this statement. 

Statement (5): from statement (3), for any first-order flex $\boldsymbol{\rho}'$ there is a self-stress $\boldsymbol{\omega}_s(\boldsymbol{\rho}')$ such that  $\boldsymbol{\rho}'^T [\boldsymbol{\omega}_s \cdot \mathrm{d}^2 \boldsymbol{A}/\mathrm{d} \boldsymbol{\rho}^2]\, \boldsymbol{\rho}' > 0$. If $\mathrm{rank}(\mathrm{d} \boldsymbol{A}/\mathrm{d} \boldsymbol{\rho})=N_c-1$, the dimension of first-order flex is 1 and the nullspace of  $\mathrm{d} \boldsymbol{A}/\mathrm{d} \boldsymbol{\rho}$ is $c\boldsymbol{\rho}_1', c \in \mathbb{R}$. The self-stress $\boldsymbol{\omega}_s(\boldsymbol{\rho}_1')$ will stabilize this rigid origami since $\boldsymbol{\rho}'^T [\boldsymbol{\omega}_s(\boldsymbol{\rho}_1') \cdot \mathrm{d}^2 \boldsymbol{A}/\mathrm{d} \boldsymbol{\rho}^2]\,\boldsymbol{\rho}'= c^2 \boldsymbol{\rho}_1'^T [\boldsymbol{\omega}_s(\boldsymbol{\rho}_1') \cdot \mathrm{d}^2 \boldsymbol{A}/\mathrm{d} \boldsymbol{\rho}^2]\, \boldsymbol{\rho}_1' > 0$. Next, if the dimension of self-stress is $1$. Denote this base as $\boldsymbol{\omega}_1$. If this rigid origami is not prestress stable, there will exist a first-order flex $\boldsymbol{\rho}'$ such that for all choice of $c$, $c \boldsymbol{\rho}'^T [\boldsymbol{\omega}_1 \cdot \mathrm{d}^2 \boldsymbol{A}/\mathrm{d} \boldsymbol{\rho}^2]\, \boldsymbol{\rho}' \le 0$. First, $\boldsymbol{\rho}'^T [\boldsymbol{\omega}_1 \cdot \mathrm{d}^2 \boldsymbol{A}/\mathrm{d} \boldsymbol{\rho}^2]\, \boldsymbol{\rho}' \neq 0$ since the rigid origami is second-order rigid. Second, by choosing $c=\pm 1$, $c \boldsymbol{\rho}'^T [\boldsymbol{\omega}_1 \cdot \mathrm{d}^2 \boldsymbol{A}/\mathrm{d} \boldsymbol{\rho}^2]\, \boldsymbol{\rho}'$ could be greater than 0. These lead to a contradiction.  
\end{proof}

\begin{rem}
Prestress stability requires a single self-stress $\boldsymbol{\omega}$ such that the quadratic form is positive for every first-order flex, while the second-order rigidity requires a ``suitable'' self-stress $\boldsymbol{\omega}(\boldsymbol{\rho}')$ for every first-order flex such that the quadratic form is positive. Physically, such a $\boldsymbol{\omega}(\boldsymbol{\rho}')$ ``blocks'' a possible second-order flex for a given first-order flex.
\end{rem}

\begin{rem}
For a planar rigid origami, after defining a reciprocal diagram, the first and second order rigid-foldability could be graphically explained as the existence, and the zero-area property, of the reciprocal diagram \cite{watanabe_method_2009,tachi_design_2012,demaine_zero-area_2016}.
\end{rem}

An attempt to find a rigid origami that is second-order rigid but not pre-stress stable is provided in Section S3.4 of the supplementary material. Here we conjecture that (\textit{regular} means the rigidity matrix has maximum rank)

\begin{conj} \label{conj: single-hole prestress stability}
A rigid but not regular single-hole is prestress stable.
\end{conj}

\section{Relation among different levels of rigidity} \label{section: relation}

In this section we will prove the relation among the rigidity discussed in the above sections, and which is illustrated in figure \ref{fig: introduction}.

\begin{thm} \label{thm: rigidity}
The relation among first-order or static rigidity, prestress stability and second-order rigidity.
\begin{enumerate} [(1)]
	\item A rigid origami is prestress stable if first-order rigid or statically rigid.
	\item A rigid origami is second-order rigid if prestress stable.
	\item A rigid origami is rigid if second-order rigid.
\end{enumerate}
\end{thm}

\begin{proof}~

Statement (1): Set $\boldsymbol{\omega}_s=\boldsymbol{0}$, the total stiffness $K=(\mathrm{d} \boldsymbol{A}/\mathrm{d} \boldsymbol{\rho})^T \cdot B \cdot \mathrm{d} \boldsymbol{A}/\mathrm{d} \boldsymbol{\rho}$ is now positive-definite.

Statement (2):  If a rigid origami is prestress stable, for any first-order flex $\boldsymbol{\rho}'$, there is a uniform $\boldsymbol{\omega}_s$ such that $\boldsymbol{\rho}'^T [\boldsymbol{\omega}_s \cdot \mathrm{d}^2 \boldsymbol{A}/\mathrm{d} \boldsymbol{\rho}^2]\, \boldsymbol{\rho}' > 0$. From statement (3) of proposition \ref{prop: 6}, this rigid origami is second-order rigid.

Statement (3): We need to prove that for rigid origami, a continuous flex implies a second-order flex. This could be done by transferring the consistency constraint $\boldsymbol{A}(\boldsymbol{\rho})=\boldsymbol{0}$ to a polynomial system $\boldsymbol{A}(\boldsymbol{t})=\boldsymbol{0}$ with the normalized folding angle expression $\boldsymbol{t}=\tan{(\boldsymbol{\rho}/2)}$, and we claim that the definitions on local rigidity are equivalent for these two expressions (details are provided in Section S5 of the supplementary material). It turns out that a continuous flex is equivalent to an analytical flex in the normalized folding angle expression \cite{he_rigid_2019}. Denote an analytical flex for a foldable rigid origami starting from $(\boldsymbol{t},\lambda')$ by $\boldsymbol{\gamma}:[0,1]\ni s \rightarrow \{\boldsymbol{t}\}$. This flex could be parametrized by a single $s$:
\begin{equation}
	\boldsymbol{\gamma}=\boldsymbol{t}+ \sum \limits_{n=1}^{\infty} \dfrac{\boldsymbol{a}_n}{n!}  s^n, \quad \boldsymbol{A}(\boldsymbol{\gamma}) \equiv \boldsymbol{0}
\end{equation}
where not all $\boldsymbol{a}_n=\boldsymbol{0}$. 

If $\boldsymbol{a}_1 \neq \boldsymbol{0}$, as we know that
\begin{equation}
	\left.\dfrac{\mathrm{d} \boldsymbol{A}}{\mathrm{d} s}\right|_{s=0}= \boldsymbol{0}, \quad \left.\dfrac{\mathrm{d}^2 \boldsymbol{A}}{\mathrm{d} s^2}\right|_{s=0}=\boldsymbol{0}
\end{equation}
then $(\boldsymbol{a}_1, \boldsymbol{a}_2)$ would be a second-order flex satisfying equation (S5.10). 

If $\boldsymbol{a}_1 = \boldsymbol{0}$, then as $\boldsymbol{\gamma} \neq \boldsymbol{0}$, there must be some first non-zero term $\boldsymbol{a}_k$.  Then, as we know that 
\begin{equation}
	\left.\dfrac{\mathrm{d}^i \boldsymbol{A}}{\mathrm{d} s ^i}\right|_{s=0}= \boldsymbol{0}, \quad 1 \le i \le 2k
\end{equation}
hence $\boldsymbol{a}_k$ would be a first-order flex, and $(\boldsymbol{a}_k, 2\boldsymbol{a}_{2k}/\tbinom{2k}{k})$ would be a second-order flex satisfying equation (S5.10). Here $\tbinom{2k}{k}$ is the binomial coefficient in combination.
\end{proof}

Note that none of the statement in theorem \ref{thm: rigidity} is reversible. Apart from the examples given in Sections S3.3 and S3.4 that are prestress stable but not first-order rigid, in Section S3.5 we show an example that is rigid but not second-order rigid.

\section{From local rigid-foldability to rigid-foldability} \label{section: extension}

When studying the hierarchical relation described in theorem \ref{thm: rigidity}, it turns out that for some rigid origami, different levels of rigidity might be equivalent. In particular, a first-order flex which will not lead to crossing of panels might be extended to a continuous flex. 

\begin{prop} \label{prop: extension} Extension of local rigid-foldability for some special rigid origami. Recall that a rigid origami is \textit{regular} if the rigidity matrix has maximum rank.
\begin{enumerate} [(1)]
	\item A regular rigid origami is rigid-foldable if first-order rigid-foldable. 
	\item A single-vertex is rigid-foldable if not prestress stable. 
	\item A planar quadrilateral mesh where each vertex is flat-foldable is at least second-order rigid-foldable.
\end{enumerate}
\end{prop}	

\begin{proof}
Statement (1): Here the rigid origami $(\boldsymbol{\rho}, \lambda')$ is also called \textit{independent}, where the only self-stress is $\boldsymbol{0}$. Suppose there are $N_v$ inner vertices, $N_h$ holes and $N_c$ inner creases. From the Implicit Function Theorem \cite[Section 8.5]{zorich_mathematical_2004}, at a neighbourhood of $(\boldsymbol{\rho}, \lambda')$, the folding angle space would be a manifold with dimension $N_c-\mathrm{rank}(\mathrm{d} \boldsymbol{A}/\mathrm{d} \boldsymbol{\rho})=N_c-3N_v-6N_h>0$, hence there would be a continuous flex starting from $(\boldsymbol{\rho}, \lambda')$.  

Statement (2): For a non-planar single-vertex, it is regular, hence either first-order rigid or rigid-foldable. For a planar single-vertex, from proposition \ref{prop: single-vertex prestress stability}, if not prestress stable it would be rigid-foldable. 

Statement (3):  If each vertex of a planar quadrilateral mesh is flat-foldable, the relations among the tangent of half of all the folding angles, shown in the consistency constraint, are linear, even though this quadrilateral mesh might not be rigid-foldable \cite{tachi_geometric_2010}. Consider the normalized folding angle expression, the consistency constraint $\boldsymbol{A}(\boldsymbol{t})=\boldsymbol{0}$ could be rewritten as a linear system $\boldsymbol{A}'(\boldsymbol{t})=\boldsymbol{0}$ among $\boldsymbol{t}$, hence the Hessian $\mathrm{d}^2 \boldsymbol{A}'/\mathrm{d} \boldsymbol{t}^2$ is zero. Since a planar quadrilateral mesh is also first-order rigid-foldable, it cannot be pre-stress stable or second-order rigid. 
\end{proof}

A rigid but not second-order rigid example for the statement (3) of proposition \ref{prop: extension} is provided in Subsection S3.5 of the supplementary material, but the stress matrix $\boldsymbol{\omega}_s \cdot \mathrm{d}^2 \boldsymbol{A}/\mathrm{d} \boldsymbol{\rho}^2$ is not zero. An explanation is, although $\boldsymbol{A}(\boldsymbol{t})=\boldsymbol{0}$ is essentially a linear system with the special choice of sector angles such that every vertex is flat-foldable, it would be in the form of a complicated polynomial system consisting the square of linear relations, hence $\boldsymbol{\omega}_s \cdot \mathrm{d}^2 \boldsymbol{A}/\mathrm{d} \boldsymbol{\rho}^2$ is not zero. Further, not every first-order flex can be extended to a second-order flex in this example when choosing the consistency constraint to be $\boldsymbol{A}(\boldsymbol{\rho})=\boldsymbol{0}$. We claim that the conclusion on local rigidity should be invariant to the choice of form of consistency constraint. 

Proposition \ref{prop: extension} opens a promising topic for a rigid origami, that is, to explore the level of local rigidity and find whether some of these levels are in fact equivalent.

As stated in the introduction, when some folding angles are $\pm \pi$, a first-order flex calculated from the independent consistency constraint $\boldsymbol{A}(\boldsymbol{\rho})=\boldsymbol{0}$ is extendable to a flex only when this first-order flex is pointing away from $\pm \pi$. Some examples on this topic are provided in Section S3.6 of the supplementary material. 

\section{Numerical methods for rigidity analysis}

In this section we will consider how to analyse the local rigidity for a rigid origami using numerical method when the size of $\mathrm{d} \boldsymbol{A}/\mathrm{d} \boldsymbol{\rho}$ is large. Several important questions are:

\begin{enumerate} [(1)]
\item How to determine the first-order rigidity of a rigid origami or find the collection of first-order flex and self-stress?
\item Will a given self-stress stabilize a rigid origami?
\item How to find the collection of self-stress that can stabilize a rigid origami?
\item Can a first-order flex be extended to a second-order flex?
\item How to find the collection of first-order flex that can be extended to a second-order flex?
\end{enumerate}

For question (1), if we know the position of each vertex, the direction vector of each inner crease can be calculated directly, and the rigidity matrix $\mathrm{d} \boldsymbol{A}/\mathrm{d} \boldsymbol{\rho}$ in a global coordinate system can be obtained by assigning entries for a sparse matrix. The next step is applying the singular value decomposition (SVD) to $\mathrm{d} \boldsymbol{A}/\mathrm{d} \boldsymbol{\rho}$ to find the information of rank, nullspace and left nullspace \cite{pellegrino_structural_1993, chan_improved_1982}.

Question (2) is the forward problem of determining prestress stability. The Hessian $\mathrm{d}^2 \boldsymbol{A}/\mathrm{d} \boldsymbol{\rho}^2$ of a rigid origami can also be obtained by assigning entries for a sparse matrix if the position of each vertex is known. With the information of the nullspace of $\mathrm{d} \boldsymbol{A}/\mathrm{d} \boldsymbol{\rho}$ calculated in (1), we need to know the positive definiteness of $[\rho']^T \cdot \boldsymbol{\omega}_s \cdot \mathrm{d}^2 \boldsymbol{A}/\mathrm{d} \boldsymbol{\rho}^2 \cdot [\rho']$ from statement (3) of proposition \ref{prop: 5}, which is symmetric. The eigenvalues and eigenvectors for a real, sparse and symmetric matrix could be found, for instance, by the modified Lanczos algorithm \cite{ruhe_implementation_1979, edwards_use_1979}.

Question (3) is the inverse problem of determining prestress stability. The collection of self-stress that can stabilize a rigid origami turns out to be an elliptic system. Suppose the bases of self-stress are $\boldsymbol{\omega}_s^i$ ($1 \le i \le 3N_v+6N_h-\mathrm{rank}(\mathrm{d} \boldsymbol{A}/\mathrm{d} \boldsymbol{\rho})$), which is calculated in (1). Now the problem becomes, is there a linear combination of these $3N_v+6N_h-\mathrm{rank}(\mathrm{d} \boldsymbol{A}/\mathrm{d} \boldsymbol{\rho})$ real, sparse and symmetric matrices $[\rho']^T \cdot \boldsymbol{\omega}_s^i \cdot \mathrm{d}^2 \boldsymbol{A}/\mathrm{d} \boldsymbol{\rho}^2 \cdot [\rho']$ that is positive-definite. This is a problem in semi-definite programming which has been well studied \cite{vandenberghe_semidefinite_1996}. We could set this problem as 
\begin{equation*}
\begin{gathered}
	\mathrm{minimize} \quad \boldsymbol{d}^T\boldsymbol{c} \\
	\mathrm{s.t.} ~ \sum c_i[\rho']^T \boldsymbol{\omega}_s^i \dfrac{\mathrm{d}^2 \boldsymbol{A}}{\mathrm{d} \boldsymbol{\rho}^2} [\rho'] \quad \mathrm{positive ~ definite}
\end{gathered}
\end{equation*}
where $\boldsymbol{d} \in \mathbb{R}^n$ is a given vector that converges the solution set of $\boldsymbol{c}$ to be elliptic. Note that even a stabilizing self-stress $\boldsymbol{\omega}_s$ is found, the proof of existence of $k$ is not constructive in proposition \ref{prop: 5}. Other techniques need to be applied to determine how small $k$ has to be. 

Question (4) is the forward problem of determining second-order rigidity. From statement (2) in proposition \ref{prop: 6}, the problem is to consider whether the given first-order flex $\boldsymbol{\rho}$' is in the quadratic nullspace of every $\boldsymbol{\omega}_s^i \cdot \mathrm{d}^2 \boldsymbol{A}/\mathrm{d} \boldsymbol{\rho}^2$, which is calculated in (3). If not, $\boldsymbol{\rho}'$ could be extended to a second-order flex. 

Question (5) is the inverse problem of determining second-order rigidity. The collection of first-order flex that can be extended to a second-order flex is also an elliptic system. From statement (4) in proposition \ref{prop: 6}, we need to find the common root of the quadratic form for each $[\rho']^T \cdot \boldsymbol{\omega}_s^i \cdot \mathrm{d}^2 \boldsymbol{A}/\mathrm{d} \boldsymbol{\rho}^2 \cdot [\rho']$, where $\{\boldsymbol{\omega}_s^i\}$ is a base of self-stress $(1 \le i \le 3N_v+6N_h-\mathrm{rank}(\mathrm{d} \boldsymbol{A}/\mathrm{d} \boldsymbol{\rho}))$. Since each $[\rho']^T \cdot \boldsymbol{\omega}_s^i \cdot \mathrm{d}^2 \boldsymbol{A}/\mathrm{d} \boldsymbol{\rho}^2 \cdot [\rho']$ is real and symmetric, we could write its eigenvalues as $s_j$ and its orthonormal vectors as $\boldsymbol{v}_j, (1 \le j \le N_c-\mathrm{rank}(\mathrm{d} \boldsymbol{A}/\mathrm{d} \boldsymbol{\rho}))$. If $\boldsymbol{a}$ is a root of the quadratic form,
\begin{equation}
\begin{gathered}
	\boldsymbol{a}=\sum_1^{N_c-\mathrm{rank}(\mathrm{d} \boldsymbol{A}/\mathrm{d} \boldsymbol{\rho})} c_j\boldsymbol{v}_j \\
	\sum_1^{N_c-\mathrm{rank}(\mathrm{d} \boldsymbol{A}/\mathrm{d} \boldsymbol{\rho})} c_j^2 s_j=0
\end{gathered}
\end{equation}
That is to say, the square of coefficients $\boldsymbol{c}^2$ when $\boldsymbol{a}=\sum c_j \boldsymbol{v}_j$ should be orthogonal to the eigenvalues $\boldsymbol{s}$. The next step is to find the intersection of such $3N_v+6N_h-\mathrm{rank}(\mathrm{d} \boldsymbol{A}/\mathrm{d} \boldsymbol{\rho})$ quadratic nullspace for each base of self-stress. For a large rigid origami the computation would be expensive.

\section{Conclusion}

We have shown that rigid origami can, with advantage, be analysed from a rigidity point  of view.  This is an inversion of the usual focus on folding.  Rather than consider when a paper can fold, we have examined various ways in which the design of an origami might prevent folding.  We think this perspective will prove to be of further use in the development of novel folding patterns, or indeed in the design of structures formed from origami where some rigidity is required.

\section*{Acknowledgement}

ZH was partially funded by a fellowship awarded by the George and Lilian Schiff Foundation. We thank Tomohiro Tachi for enlightening discussions in a OIST workshop on origami and deployable mechanisms in May 2019.

\vskip2pt

\bibliographystyle{unsrt}

\bibliography{Rigid-Folding}

\end{document}


\maketitle
	
\numberwithin{equation}{section}
\numberwithin{prop}{section}

In this supplementary material we will provide the preliminaries for a bar-joint framework in Section \ref{section: s1}, then give the detailed derivation for the first-order and second-order derivative of the independent consistency constraint $\boldsymbol{A}(\boldsymbol{\rho})=\boldsymbol{0}$ in Section \ref{section: s2}. Some examples on writing the rigidity matrix for a large rigid origami; calculating internal force and self-stress; determining prestress stability and second-order rigidity; and special cases when some folding angles are $\pm \pi$, are listed in Section \ref{section: s3}. The proof of prestress stability for a rigid planar single-vertex is given in Section \ref{section: example prestress vertex}. The definitions and main conclusions for local rigidity in the normalized folding angle expression $\boldsymbol{t}=\tan (\boldsymbol{\rho}/2)$ are provided in Section \ref{section: normalized}.

\section{Preliminaries on a bar-joint framework} \label{section: s1}

Here we will introduce the basics of rigidity theory for a bar-joint framework, as it may prove useful background for a reader new to rigidity theory.

\begin{defn}
	A simple graph $G$ can be written as a pair of its vertices and edges $G=(V,E)$. Suppose there are $v$ vertices and $e$ edges. A \textit{configuration} $\boldsymbol{p}=(p_1, p_2, ... , p_v)$ is the position of all vertices, such that $p_i \neq p_j$ for each $\{i,j\} \in E$. A \textit{bar-joint framework} $G(\boldsymbol{p})$ is a graph $G$ and a \textit{configuration} $\boldsymbol{p}$. Here a vertex is called a \textit{joint}, and an edge is called a \textit{bar}. Given a graph $G$ and the length of all edges $l_{ij}$, if there is a configuration $\boldsymbol{p}$, such that $||p_i-p_j||=l_{ij}$ for each $\{i,j\} \in E$, this $\boldsymbol{p}$ is called a \textit{realization} of $G$. If there exists a deleted neighbourhood of a realization $\boldsymbol{p}$ such that no configuration is a realization, then $\boldsymbol{p}$ is \textit{rigid}, otherwise $\boldsymbol{p}$ is \textit{flexible}.
\end{defn}

The geometry constraint of a bar-joint framework is $||p_i-p_j||=l_{ij}$, and it is natural to calculate the Jacobian for the information of first-order rigidity. The Jacobian, i.e, \textit{rigidity matrix} $R_G$ of a realization $\boldsymbol{p}$ with graph $G$ is the matrix form of:
\begin{equation}
	(p_i-p_j)^T(p_i'-p_j')=0
\end{equation}

\begin{defn}
	At a realization $\boldsymbol{p}$ of a given bar-joint framework, a \textit{first-order flex} $\boldsymbol{p}'$ is a vector in the tangent space, which can be expressed as $R_G \cdot \boldsymbol{p}'=\boldsymbol{0}$. A \textit{trivial first-order flex} is the ``speed'' of rotation and translation among all the first-order flexes, which can be written as $\forall i, p_i'=Sp_i+t$. Here $S$ is a constant skew-symmetric matrix (for rotation) and $t$ is a constant vector (for translation). If the solution space of $R_G \cdot \boldsymbol{p}'=\boldsymbol{0}$ is the collection of trivial first-order flexes, or equivalently, $\mathrm{rank} (R_G)=3v-6$, the framework is said to be \textit{first-order rigid} at $\boldsymbol{p}$, otherwise \textit{first-order flexible}.
\end{defn}

Next we suppose all the loads are forces concentrated on joints. Consider the equilibrium of each joint, we have the following on the static rigidity of a bar-joint framework.

\begin{defn}
	A bar-joint framework can \textit{resolve} a set of equilibrium forces applied on its joints at a realization $\boldsymbol{p}$ if each joint is in equilibrium, i.e. the following equation has a solution
	\begin{equation}
		\sum \limits_j \sigma_{ij}(p_i-p_j)=F_i
	\end{equation}
	where $\sigma_{ij}$ is the axial force per length (or called \textit{stress}) on the bar joining vertices $i$ and $j$, $F_i$ is the external force applied on vertex $i$. If we denote $\boldsymbol{\omega}=\{\omega_{ij}\}$ (in total $e$ entries) and $\boldsymbol{F}=\{F_i\}$, the matrix form of the above equation is:
	\begin{equation}
		\boldsymbol{\omega} R_G=\boldsymbol{F}
	\end{equation}
	A bar-joint framework is \textit{statically rigid} at a realization $\boldsymbol{p}$ if it can resolve every set of equilibrium forces applied on its joints. When $\boldsymbol{F}=\boldsymbol{0}$, $\boldsymbol{\omega}$ is called the \textit{self-stress}, denoted by $\boldsymbol{\omega}_s$.
	\begin{equation}
		\boldsymbol{\omega}_s R_G=\boldsymbol{0}
	\end{equation}
\end{defn}

\begin{prop}
	The following statements for a realization $\boldsymbol{p}$ of a bar-joint framework $G$ are equivalent:
	\begin{enumerate} [(1)]
		\item $G(\boldsymbol{p})$ is first-order rigid.
		\item $G(\boldsymbol{p})$ is statically rigid.
		\item The dimension of the collection of self-stress at $\boldsymbol{p}$ is $e-3v+6$.
	\end{enumerate}
\end{prop}

\begin{proof}
	For the rigidity matrix $R_G$, a zero nullspace is equivalent to a full image. The rank of its left nullspace is $e-3v+6$ when the nullspace is zero.
\end{proof}

We might also be interested in the behaviour of a bar-joint framework under a general form of load, which means, a random set of equilibrium forces and torques (either discrete or continuous). For each bar, the load applied could be decomposed into the load applied on its joints, hence for a general load the behaviour could be obtained. Further, it is also possible to introduce axial, shear, bending and torsional stiffness at each bar, and possibly built-in joint instead of simply supported joint for some bars, then study its behaviour under a general load. This is the ``displacement method'' introduced in structures engineering textbook for small deformation. An alternative way is to write the elastic potential energy as the function of joint displacement and rotation; and the external potential energy as a functional of deflection and torsion of each bar, then apply the variation method to find the equilibrium and stability condition with a set of chosen boundary condition. 

\section{First-order and second-order derivative} \label{section: s2}

\begin{figure}[!tb]
	\noindent \begin{centering}
		\includegraphics[width=1\linewidth]{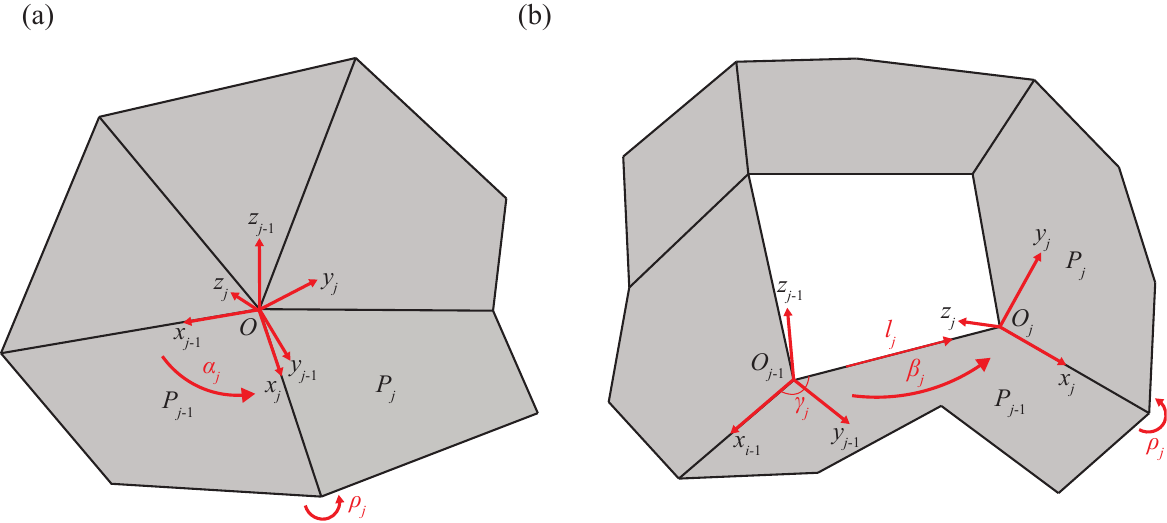}
		\par\end{centering}
	
	\caption{\label{fig: single vertex and hole}(a) and (b) show the rotation and possible translation of local coordinate systems around a degree-5 vertex and a degree-5 hole. The first Betti number of (a) is 0 and of (b) is 1. The form of consistency constraint is with equation (2.1) for (a) and with equation (2.2) for (b) in the main text.}
\end{figure}

This section returns to rigid origami, and gives the derivations for the first-order and second-order derivative of the independent consistency constraint $\boldsymbol{A}(\boldsymbol{\rho})=\boldsymbol{0}$ for a single-vertex or hole.

First, a series of local coordinate systems are attached to each panel of a single-vertex or single hole (figure \ref{fig: single vertex and hole}). Given a global coordinate system for a rigid origami, consider its restriction to a degree-$n$ single-vertex, $\boldsymbol{x}_j=[x_{1j};x_{2j};x_{3j}]$, $\boldsymbol{y}_j=[y_{1j};y_{2j};y_{3j}]$, $\boldsymbol{z}_j=[z_{1j};z_{2j};z_{3j}]$ are the direction vectors of axes $x_j,y_j,z_j$ of the local coordinate system attached to panel $j (1 \le j \le n)$ measured in the global coordinate system. We further write them collectively as $X_j$,
\begin{equation}
	X_j=[\boldsymbol{x}_j,\boldsymbol{y}_j,\boldsymbol{z}_j]=\left[ \begin{array}{ccc}
		x_{1j} & y_{1j} & z_{1j}\\
		x_{2j} & y_{2j} & z_{2j}\\
		x_{3j} & y_{3j} & z_{3j}
	\end{array} \right]
\end{equation}
At a rigidly folded state $\boldsymbol{\rho}$, $R$ in equation (2.1) satisfies:
\begin{equation} \label{eq: 22}
	R=\prod_{i=1}^{n}R_i=I=X_nX_n^T=X_nRX_n^T
\end{equation}
All the products used here are post-multiplication, then ($\prod \limits_{i=j}^{1}$ means product from subscript $j$ to $1$)
\begin{equation} \label{eq: Jacobian 3}
	\begin{gathered}
		\dfrac{\partial R}{\partial \rho_j}=X_n\dfrac{\partial R}{\partial \rho_j}X_n^T=X_n\left(\prod_{i=1}^{j-1}R_i\right)\dfrac{\mathrm{d} R_j}{\mathrm{d} \rho_j}\left(\prod_{i=j+1}^{n}R_i\right)X_n^T \\
		=X_n\left(\prod_{i=1}^{j}R_i\right)R_j^T\dfrac{\mathrm{d} R_j}{\mathrm{d} \rho_j}\left(\prod_{i=j}^{1}R_i^T\right) X_n^T=\left(X_n\prod_{i=1}^{j}R_i\right)
		\left[ \begin{array}{ccc}
			0 & 0 & 0\\
			0 & 0 & -1\\
			0 & 1 & 0
		\end{array} \right]
		\left(X_n\prod_{i=1}^{j}R_i\right)^T \\
		=X_j
		\left[ \begin{array}{ccc}
			0 & 0 & 0\\
			0 &  0 & -1\\
			0                &                 1 & 0
		\end{array} \right]
		X_j^T =
		\left[ \begin{array}{cccc}
			0 & -x_{3j} & x_{2j} \\
			x_{3j} & 0 & -x_{1j} \\
			-x_{2j} & x_{1j} & 0 
		\end{array} \right]
	\end{gathered}
\end{equation}
This leads to the first-order derivative of independent consistency constraint $\boldsymbol{A}(\boldsymbol{\rho})=\boldsymbol{0}$ for a single-vertex in equation (3.1).

Next, $T$ in equation (2.2) can be written as
\begin{align} 
	T=\prod_{i=1}^{n}
	\left[ \begin{array}{cccc}
		R_i & \boldsymbol{v}_i \\
		0 & 1 
	\end{array} \right]=I
\end{align}
where $\boldsymbol{v}_i=[l_i \cos \gamma_i;l_i \sin \gamma_i;0]$ is the position of $O_i$ measured in the local coordinate system built on panel $i-1$, therefore
\begin{equation} \label{eq: 2}
	\prod_{i=1}^{n}R_i=I
\end{equation}
\begin{equation} \label{eq: 1}
	\sum_{i=1}^{n} \prod_{k=1}^{i-1}R_kv_i=0=X_n\sum_{i=1}^{n} \prod_{k=1}^{i-1}R_kv_i
\end{equation}
The derivative of the left hand side of equation \eqref{eq: 1} with respect to $\rho_j$ is
\begin{equation} \label{eq: 33}
	\begin{gathered}
		\dfrac{\partial}{\partial \rho_j}\left(\sum_{i=1}^{n} \prod_{k=1}^{i-1}R_kv_i\right)=\dfrac{\partial}{\partial \rho_j}\left(X_n\sum_{i=1}^{n} \prod_{k=1}^{i-1}R_kv_i\right) \\
		=\sum_{i=j+1}^{n} \left(X_n\prod_{k=1}^{j-1} R_k\right) \dfrac{\mathrm{d} R_j}{\mathrm{d} \rho_j} \prod_{k=j+1}^{i-1}R_k v_i = \left(X_n\prod_{k=1}^{j-1} R_k\right) \dfrac{\mathrm{d} R_j}{\mathrm{d} \rho_j} \sum_{i=j+1}^{n} \prod_{k=j+1}^{i-1}R_k v_i
	\end{gathered}
\end{equation}
We also have
\begin{equation} \label{eq: 44}
	\begin{gathered}
		\sum_{i=j+1}^{n} \prod_{k=j+1}^{i-1}R_k v_i= \left(\prod_{k=j}^{1}R_k^T \right) \sum_{i=j+1}^{n} \prod_{k=1}^{i-1}R_k v_i=- \left(\prod_{k=j}^{1}R_k^T \right) \sum_{i=1}^{j} \prod_{k=1}^{i-1}R_kv_i
	\end{gathered}
\end{equation}
hence 
\begin{equation}
	\begin{gathered}
		\dfrac{\partial}{\partial \rho_j}\left(\sum_{i=1}^{n} \prod_{k=1}^{i-1}R_kv_i\right)=-\left(X_n\prod_{k=1}^{j-1}R_k\right)\dfrac{\mathrm{d} R_j}{\mathrm{d} \rho_j}\left(\prod_{k=j}^{1}R_k^T\right)\sum_{i=1}^{j} \prod_{k=1}^{i-1}R_kv_i \\
		=-\left(X_n\prod_{k=1}^{j-1}R_k\right)\dfrac{\mathrm{d} R_j}{\mathrm{d} \rho_j}\left(\prod_{k=j}^{1}R_k^TX_n^T\right)\left(X_n\sum_{i=1}^{j} \prod_{k=1}^{i-1}R_kv_i\right)=-\dfrac{\partial R}{\partial \rho_j} \left(X_n\sum_{i=1}^{j} \prod_{k=1}^{i-1}R_kv_i\right) \\
		=-\left[ \begin{array}{cccc}
			0 & -x_{3j} & x_{2j} \\
			x_{3j} & 0 & -x_{1j} \\
			-x_{2j} & x_{1j} & 0 
		\end{array} \right] \left(X_n\sum_{i=1}^{j} \prod_{k=1}^{i-1}R_kv_i\right) = \boldsymbol{O}_j \times \boldsymbol{x}_j
	\end{gathered}
\end{equation}
where 
\begin{equation}
	\boldsymbol{O}_j=X_n\sum_{i=1}^{j} \prod_{k=1}^{i-1}R_kv_i
\end{equation} 
is the position of origin $O_j$ measured in the global coordinate system. This leads to the first-order derivative of independent consistency constraint $\boldsymbol{A}(\boldsymbol{\rho})=\boldsymbol{0}$ for a single-hole in equation (3.2). 

We now derive the Hessian $\mathrm{d}^2 \boldsymbol{A}/\mathrm{d} \boldsymbol{\rho}^2$ for a single-vertex and hole. Recall equations \eqref{eq: 22} and \eqref{eq: Jacobian 3}, for $1 \le k < j \le n$,
\begin{equation} \label{eq: 5}
	\begin{gathered}
		\dfrac{\partial^2 R}{\partial \rho_k \partial \rho_j}=X_n\dfrac{\partial^2 R}{\partial \rho_k \partial \rho_j}X_n^T=X_n\left(\prod_{i=1}^{k-1}R_i\right) \dfrac{\mathrm{d} R_k}{\mathrm{d} \rho_k} \left(\prod_{i=k+1}^{j-1}R_i\right) \dfrac{\mathrm{d} R_j}{\mathrm{d} \rho_j} \left(\prod_{i=j+1}^{n}R_i \right) X_n^T
	\end{gathered}
\end{equation}
We also have
\begin{equation} 
	\begin{gathered}
		\prod_{i=k+1}^{j-1}R_i=\prod_{i=k}^{1}R_i^T\prod_{i=n}^{j}R_i^T=\prod_{i=k}^{1}R_i^T\prod_{i=1}^{j-1}R_i
	\end{gathered}
\end{equation}
then
\begin{equation} \label{eq: 55}
	\begin{gathered}
		\dfrac{\partial^2 R}{\partial \rho_k \partial \rho_j}=X_n\left(\prod_{i=1}^{k}R_i\right) R_k^T\dfrac{\mathrm{d} R_k}{\mathrm{d} \rho_k} \left(\prod_{i=k}^{1}R_i^T \right)X_n^TX_n \left(\prod_{i=1}^{j}R_i\right)R_j^T\dfrac{\mathrm{d} R_j}{\mathrm{d} \rho_j}\left(\prod_{i=j}^{1}R_i^T\right)X_n^T \\
		=\left(X_n\prod_{i=1}^{k}R_i\right) \left[ \begin{array}{ccc}
			0 & 0 & 0\\
			0 & 0 & -1\\
			0 & 1 & 0
		\end{array} \right] \left(X_n\prod_{i=1}^{k}R_i\right)^T  \left(X_n\prod_{i=1}^{j}R_i\right)\left[ \begin{array}{ccc}
			0 & 0 & 0\\
			0 & 0 & -1\\
			0 & 1 & 0
		\end{array} \right]\left(X_n\prod_{i=1}^{j}R_i\right)^T \\
		=\left[ \begin{array}{cccc}
			0 & -x_{3k} & x_{2k} \\
			x_{3k} & 0 & -x_{1k} \\
			-x_{2k} & x_{1k} & 0 
		\end{array} \right]\left[ \begin{array}{cccc}
			0 & -x_{3j} & x_{2j} \\
			x_{3j} & 0 & -x_{1j} \\
			-x_{2j} & x_{1j} & 0 
		\end{array} \right] \\
		=\left[ \begin{array}{cccc}
			-x_{2k}x_{2j}-x_{3k}x_{3j} & x_{2k}x_{1j} & x_{3k}x_{1j} \\
			x_{1k}x_{2j} & -x_{3k}x_{3j}-x_{1k}x_{1j} & x_{3k}x_{2j} \\
			x_{1k}x_{3j} & x_{2k}x_{3j} & -x_{1k}x_{1j}-x_{2k}x_{2j} 
		\end{array} \right] 
	\end{gathered}
\end{equation}

When $k=j$, 
\begin{equation} 
	\begin{gathered}
		\dfrac{\partial^2 R}{\partial \rho_j^2}=X_n\dfrac{\partial^2 R}{\partial \rho_j^2}X_n^T=X_n\left(\prod_{i=1}^{j-1}R_i\right)\dfrac{\mathrm{d}^2 R_j}{\mathrm{d} \rho_j^2}\left(\prod_{i=j+1}^{n}R_i\right)X_n^T=X_n\left(\prod_{i=1}^{j-1}R_i\right)\dfrac{\mathrm{d}^2 R_j}{\mathrm{d} \rho_j^2}\left(\prod_{i=j}^{1}R_i^T\right)X_n^T \\
		=X_n\left(\prod_{i=1}^{j}R_i\right)R_j^T\dfrac{\mathrm{d}^2 R_j}{\mathrm{d} \rho_j^2}\left(\prod_{i=j}^{1}R_i^T\right)X_n^T =\left(X_n\prod_{i=1}^{j}R_i\right)
		\left[ \begin{array}{ccc}
			0 & 0 & 0\\
			0 & -1 & 0\\
			0 & 0 & -1
		\end{array} \right]
		\left(X_n\prod_{i=1}^{j}R_i\right)^T
		\\=\left[ \begin{array}{ccc}
			x_{1j} & y_{1j} & z_{1j}\\
			x_{2j} & y_{2j} & z_{2j}\\
			x_{3j} & y_{3j} & z_{3j}
		\end{array} \right]
		\left[ \begin{array}{ccc}
			0 & 0 & 0\\
			0 &  -1 & 0\\
			0                &                 0 & -1
		\end{array} \right]
		\left[ \begin{array}{ccc}
			x_{1j} & x_{2j} & x_{3j}\\
			y_{1j} & y_{2j} & y_{3j}\\
			z_{1j} & z_{2j} & z_{3j}
		\end{array} \right] \\ =
		\left[ \begin{array}{cccc}
			-x_{2j}^2-x_{3j}^2 & x_{2j}x_{1j} & x_{3j}x_{1j} \\
			x_{1j}x_{2j} & -x_{3j}^2-x_{1j}^2 & x_{3j}x_{2j} \\
			x_{1j}x_{3j} & x_{2j}x_{3j} & -x_{1j}^2-x_{2j}^2 
		\end{array} \right]
	\end{gathered}
\end{equation}
and when $1 \le j < k \le n$ just switch the subscripts. This leads to the second-order derivative for independent consistency constraint $\boldsymbol{A}(\boldsymbol{\rho})=\boldsymbol{0}$ for a single-vertex in equation (5.14).

Next, recall equations \eqref{eq: 2}, \eqref{eq: 1} and \eqref{eq: 33}, the second-order derivative of the left hand side of equation \eqref{eq: 1} when $1 \le k < j \le n$ is,
\begin{equation}
	\begin{gathered}
		\dfrac{\partial^2}{\partial \rho_k \partial \rho_j}\left(\sum_{i=1}^{n} \prod_{l=1}^{i-1}R_lv_i\right)=\dfrac{\partial^2}{\partial \rho_k \partial \rho_j}\left(X_n\sum_{i=1}^{n} \prod_{l=1}^{i-1}R_lv_i\right) \\
		=X_n\dfrac{\partial}{\partial \rho_k}\left[\sum_{i=j+1}^{n} \left(\prod_{l=1}^{j-1} R_l \right) \dfrac{\mathrm{d} R_j}{\mathrm{d} \rho_j} \prod_{l=j+1}^{i-1}R_l v_i\right]\\
		=X_n\sum_{i=j+1}^{n} \left(\prod_{l=1}^{k-1} R_l\right) \dfrac{\mathrm{d} R_k}{\mathrm{d} \rho_k} \left(\prod_{l=k+1}^{j-1} R_l\right) \dfrac{\mathrm{d} R_j}{\mathrm{d} \rho_j} \prod_{l=j+1}^{i-1}R_l v_i \\
		=X_n\left(\prod_{l=1}^{k-1} R_l\right) \dfrac{\mathrm{d} R_k}{\mathrm{d} \rho_k} \left(\prod_{l=k+1}^{j-1} R_l \right)\dfrac{\mathrm{d} R_j}{\mathrm{d} \rho_j} \sum_{i=j+1}^{n} \prod_{l=j+1}^{i-1}R_l v_i
	\end{gathered}
\end{equation}
apply equations \eqref{eq: 44}, \eqref{eq: 5} and \eqref{eq: 55} 
\begin{equation}
	\begin{gathered}
		\dfrac{\partial^2}{\partial \rho_k \partial \rho_j}\left(\sum_{i=1}^{n} \prod_{l=1}^{i-1}R_lv_i\right)
		=-X_n\left(\prod_{l=1}^{k-1} R_l\right) \dfrac{\mathrm{d} R_k}{\mathrm{d} \rho_k} \left(\prod_{l=k+1}^{j-1} R_l \right) \dfrac{\mathrm{d} R_j}{\mathrm{d} \rho_j} \left(\prod_{l=j}^{1}R_l^T \right) \sum_{i=1}^{j} \prod_{l=1}^{i-1}R_lv_i \\
		=-\left[ \begin{array}{cccc}
			0 & -x_{3k} & x_{2k} \\
			x_{3k} & 0 & -x_{1k} \\
			-x_{2k} & x_{1k} & 0 
		\end{array} \right]\left[ \begin{array}{cccc}
			0 & -x_{3j} & x_{2j} \\
			x_{3j} & 0 & -x_{1j} \\
			-x_{2j} & x_{1j} & 0 
		\end{array} \right] \left(X_n \sum_{i=1}^{j} \prod_{l=1}^{i-1}R_lv_i\right) = \boldsymbol{x}_k \times (\boldsymbol{O}_j \times \boldsymbol{x}_j) 
	\end{gathered}
\end{equation}

When $k=j$, the result will be $\boldsymbol{x}_j \times (\boldsymbol{O}_j \times \boldsymbol{x}_j)$, and when $1 \le j < k \le n$ just switch the subscripts. This leads to the second-order derivative for independent consistency constraint $\boldsymbol{A}(\boldsymbol{\rho})=\boldsymbol{0}$ for a single-hole in equation (5.15).

\section{Examples} \label{section: s3}

Here we give the necessary calculations for examples mentioned in the main paper.

\subsection{Rigidity matrix for a large rigid origami} \label{section: example rigidity matrix}

\begin{figure}[!tb]
	\noindent \begin{centering}
		\includegraphics[width=1\linewidth]{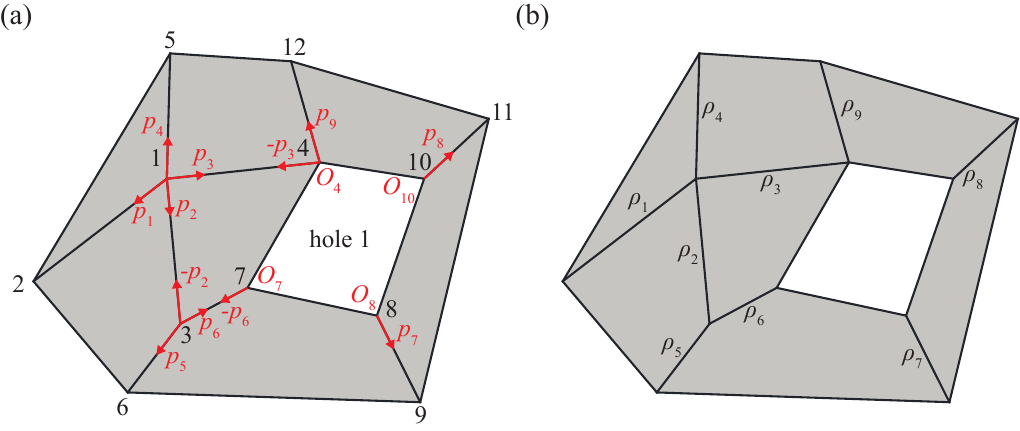}
		\par\end{centering}
	
	\caption{\label{fig: example}An example rigid origami with all the vertices labelled. The direction vectors of each inner crease $\boldsymbol{p}_j$ ($1 \le j \le 9$) and position of vertices on boundary of the hole $\boldsymbol{O}_i$ ($i=4,7,8,10$) in a global coordinate system are provided in (a). The labelling of folding angles is provided in (b).}
\end{figure}

First we provide an example of writing the rigidity matrix for a complete rigid origami in a global coordinate system, for the paper shown in figure~\ref{fig: example}. The inner creases are numbered, and each inner crease $j$ is assigned a direction vector $\boldsymbol{p}_j$.  Then we number all the vertices and write down the \textit{incidence matrix} $D$ describing the relationship between creases and vertices. If vertex $i$ is incident to inner crease $j$ and the direction vector goes out from $i$, $D_{ij}=1$; if the direction goes toward $i$, $D_{ij}=-1$, otherwise $D_{ij}=0$. $D$ is a sparse matrix.

\begin{equation} \label{eq: incidence}
	\begin{gathered}
		D=
		\left[ \scalemath{1}{{\begin{array}{cccccccccccc}
					1 & 1 & 1 & 1 & 0 & 0 & 0 & 0 & 0\\
					-1 & 0 & 0 & 0 & 0 & 0 & 0 & 0 & 0 \\
					0 & -1 & 0 & 0 & 1 & 1 & 0 & 0 & 0 \\
					0 & 0 & -1 & 0 & 0 & 0 & 0 & 0 & 1 \\
					0 & 0 & 0 & -1 & 0 & 0 & 0 & 0 & 0 \\
					0 & 0 & 0 & 0 & -1 & 0 & 0 & 0 & 0 \\
					0 & 0 & 0 & 0 & 0 & -1 & 0 & 0 & 0 \\
					0 & 0 & 0 & 0 & 0 & 0 & 1 & 0 & 0 \\
					0 & 0 & 0 & 0 & 0 & 0 & -1 & 0 & 0 \\
					0 & 0 & 0 & 0 & 0 & 0 & 0 & 1 & 0 \\
					0 & 0 & 0 & 0 & 0 & 0 & 0 & -1 & 0 \\
					0 & 0 & 0 & 0 & 0 & 0 & 0 & 0 & -1 \\
		\end{array}}} \right]
	\end{gathered}
\end{equation}

Here the inner vertices are $1$ and $3$. We take out the rows of $D$ corresponding with these inner vertices to give
\begin{equation} 
	D_{\text{vertex}}=\left[ {\begin{array}{cccccccccccc}
			1 & 1 & 1 & 1 & 0 & 0 & 0 & 0 & 0\\
			0 & -1 & 0 & 0 & 1 & 1 & 0 & 0 & 0
	\end{array}} \right]
\end{equation}
The vertices on boundary of hole $1$ (here, the only hole) are $\{4,7,8,10\}$. We sum the rows of $D$ corresponding to vertices on the boundary of each hole to give
\begin{equation} 
	D_{\text{hole}}=\left[ {\begin{array}{cccccccccccc}
			0 & 0 & -1 & 0 & 0 & -1 & 1 & 1 & 1\\
	\end{array}} \right]
\end{equation}

From equations (3.3) and (3.4), the rigidity matrix $\mathrm{d} \boldsymbol{A}/\mathrm{d} \boldsymbol{\rho}$ is obtained by replacing each $\pm 1$ by $\pm\boldsymbol{p}_j$ in column $j$ in $D_{\text{vertex}}$, and replacing $\pm 1$ by $\pm[\boldsymbol{p}_j; \boldsymbol{O}_j \times \boldsymbol{p}_j]$ in column $j$ in $D_{\text{hole}}$, then stitching these matrices together to give
\begin{equation} \label{eq: example jacobian}
	\dfrac{\mathrm{d} \boldsymbol{A}}{\mathrm{d} \boldsymbol{\rho}}=
	\left[ \scalemath{1}{{\begin{array}{cccccccccccc}
				\boldsymbol{p}_1 & \boldsymbol{p}_2 & \boldsymbol{p}_3 & \boldsymbol{p}_4 & \boldsymbol{0} & \boldsymbol{0} & \boldsymbol{0} & \boldsymbol{0} & \boldsymbol{0}\\
				\boldsymbol{0} & -\boldsymbol{p}_2 & \boldsymbol{0} & \boldsymbol{0} & \boldsymbol{p}_5 & \boldsymbol{p}_6 & \boldsymbol{0} & \boldsymbol{0} & \boldsymbol{0}\\
				\boldsymbol{0} & \boldsymbol{0} & -\boldsymbol{p}_3 & \boldsymbol{0} & \boldsymbol{0} & -\boldsymbol{p}_6 & \boldsymbol{p}_7 & \boldsymbol{p}_8 & \boldsymbol{p}_9 \\
				\boldsymbol{0} & \boldsymbol{0} & - \boldsymbol{O}_4 \times \boldsymbol{p}_3 & \boldsymbol{0} & \boldsymbol{0} & -\boldsymbol{O}_7 \times \boldsymbol{p}_6  & \boldsymbol{O}_{8} \times \boldsymbol{p}_{7} & \boldsymbol{O}_{10} \times \boldsymbol{p}_{8} & \boldsymbol{O}_4 \times \boldsymbol{p}_9 & \\
	\end{array}}} \right]
\end{equation}
For this example, with $N_v=2$ inner vertices, $N_h=1$ hole, and $N_c=9$ inner creases the size of $\mathrm{d} \boldsymbol{A}/\mathrm{d} \boldsymbol{\rho}$ is $(3N_v+6N_h) \times N_c = 12 \times 9$.

\subsection{Calculating internal force and self-stress} \label{section: example stress}

\begin{figure}[!tb]
	\noindent \begin{centering}
		\includegraphics[width=0.9\linewidth]{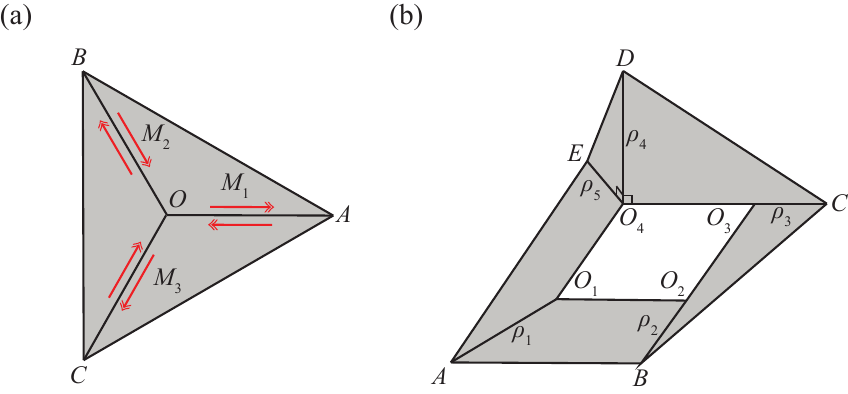}
		\par\end{centering}
	
	\caption{\label{fig: stress}(a) is a planar degree-3 single-vertex; (b) is a degree-5 single-hole. Here the load for each inner crease is a pair of opposite torques applied on adjacent panels. The internal force is a torque on each inner vertex, and a torque and a force on each hole. The equilibrium on each inner crease is shown in equations (4.6) and (4.8) in the main text.}
\end{figure}

Next we provide two simple examples showing the calculation of internal forces and states of self-stress in figure \ref{fig: stress}. Here (a) is a planar degree-3 single-vertex. The coordinates are $O(0,0,0)$, $A(1,0,0)$, $B(-1/2, \sqrt{3}/2,0)$, $C(-1/2, -\sqrt{3}/2,0)$. The incidence matrix is
\begin{equation}
	D=
	\left[
	{\begin{array}{ccc}
			1 & 1 & 1 \\
			-1 & 0 & 0 \\
			0 & -1 & 0 \\
			0 & 0 & -1 
	\end{array}} \right]
\end{equation}
The rigidity matrix is:
\begin{equation}
	\dfrac{\mathrm{d} \boldsymbol{A}}{\mathrm{d} \boldsymbol{\rho}}=
	\left[
	{\begin{array}{ccc}
			1 & -\frac{1}{2} & -\frac{1}{2} \\
			0 & \frac{\sqrt{3}}{{2}} & -\frac{\sqrt{3}}{{2}} \\
			0 & 0 & 0 
	\end{array}} \right]
\end{equation}
The direction vector of each inner crease points outside this single-vertex. Torques that are of the same magnitude but different direction are applied on a pair of panels adjacent to an inner crease. The load shown here is positive, since it does positive work on a positive increase of folding angle. Because $\mathrm{rank}(\mathrm{d} \boldsymbol{A}/\mathrm{d} \boldsymbol{\rho})=2$, here both the first-order flex and self-stress are of dimension 1, and the load that can be resolved is of dimension 2. The first-order flex is:
\begin{equation}
	\boldsymbol{\rho}'=
	a_1\left[
	{\begin{array}{c}
			1 \\
			1 \\
			1  
	\end{array}} \right]
\end{equation}
$a_1 \in \mathbb{R}$. The load that can be resolved is in the form of
\begin{equation} 
	\boldsymbol{l}=
	\left[
	{\begin{array}{c}
			l_1 \\
			l_2 \\
			-l_1-l_2  
	\end{array}} \right]
\end{equation}
$l_1, l_2 \in \mathbb{R}$, and the internal force is:
\begin{equation}
	\boldsymbol{\omega}=
	\left[
	{\begin{array}{ccc}
			l_1 & \frac{\sqrt{3}}{3}(l_1+2l_2) & b_1  
	\end{array}} \right]
\end{equation}
$b_1 \in \mathbb{R}$, then the self-stress is
\begin{equation}
	\boldsymbol{\omega}_s=
	b_1\left[
	{\begin{array}{ccc}
			0  & 0 & 1 
	\end{array}} \right]
\end{equation}

Figure \ref{fig: stress}(b) shows a degree-5 single-hole. The coordinates are $O_1(0,0,0)$, $O_2(1,0,0)$, $O_3(1,1,0)$, $O_4(0,1,0)$, $A(-\sqrt{2}/2,-\sqrt{2}/2,0)$, $B(1,-\sqrt{2}/2,0)$, $C(1+\sqrt{2}/2,0,0)$, $D(0,1,1)$, $E(-\sqrt{2}/2,1+\sqrt{2}/2,0)$. The incidence matrix is 
\begin{equation}
	D=
	\left[
	{\begin{array}{ccccc}
			1 & 0 & 0 & 0 & 0 \\
			0 & 1 & 0 & 0 & 0\\
			0 & 0 & 1 & 0 & 0\\
			0 & 0 & 0 & 1 & 1\\
			-1 & 0 & 0 & 0 & 0 \\
			0 & -1 & 0 & 0 & 0\\
			0 & 0 & -1 & 0 & 0 \\
			0 & 0 & 0 & -1 & 0\\
			0 & 0 & 0 & 0 & -1  
	\end{array}} \right]
\end{equation}
The rigidity matrix is:
\begin{equation}
	\dfrac{\mathrm{d} \boldsymbol{A}}{\mathrm{d} \boldsymbol{\rho}}=
	\left[
	{\begin{array}{ccccc}
			-\frac{\sqrt{2}}{2} & 0 & 1 & 0 & -\frac{\sqrt{2}}{2}\\
			-\frac{\sqrt{2}}{2} & -1 & 0 & 0 & \frac{\sqrt{2}}{2}\\
			0 & 0 & 0 & 1 & 0 \\
			0 & 0 & 0 & 1 & 0 \\
			0 & 0 & 0 & 0 & 0 \\
			0 & -1 & -1 & 0 & \frac{\sqrt{2}}{2}
	\end{array}} \right]
\end{equation}
Since $\mathrm{rank}(\mathrm{d} \boldsymbol{A}/\mathrm{d} \boldsymbol{\rho})=4$,  the first-order flex is of dimension 1, the self-stress is of dimension 2, and the load that can be resolved is of dimension 4. The first-order flex is:
\begin{equation}
	\boldsymbol{\rho}'=
	a_1\left[
	{\begin{array}{c}
			-\frac{\sqrt{2}}{2} \\
			\frac{1}{2} \\
			-\frac{1}{2} \\
			0 \\
			0
	\end{array}} \right]
\end{equation}
$a_1 \in \mathbb{R}$. The load that can be resolved is in the form of
\begin{equation}
	\boldsymbol{l}=
	\left[
	{\begin{array}{ccccc}
			l_1 & l_2 & -\sqrt{2}l_1+l_2 & l_3 & l_4  
	\end{array}} \right]
\end{equation}
$l_1, l_2, l_3, l_4 \in \mathbb{R}$, and the internal force is:
\begin{equation}
	\boldsymbol{\omega}=
	\left[
	{\begin{array}{cccccc}
			-l_2-\sqrt{2}l_4 & -\sqrt{2}l_1+l_2+\sqrt{2}l_4 & \sqrt{2}b_1/2 & -\sqrt{2}b_1/2+l_3 & b_2 & \sqrt{2}l_1-2l_2-\sqrt{2}l_4 \\
	\end{array}} \right]
\end{equation}
Here the coefficient for $b_1, b_2 \in \mathbb{R}$ are adjusted to have same weight, then the self-stress is:
\begin{equation}
	\boldsymbol{\omega}_s=
	\left[
	{\begin{array}{cccccc}
			0 & 0 & \sqrt{2}b_1/2 & -\sqrt{2}b_1/2 & b_2 & 0
	\end{array}} \right]
\end{equation}

\subsection{Prestress stable but not first-order rigid; the effect of load} \label{section: example stability}

We will revisit the examples in the above section to analyse their prestress stability and stability under load. Figure \ref{fig: stress}(a) shows the simplest rigid origami that is prestress stable but not first-order rigid. For this degree-3 vertex, the Hessian is:
\begin{equation}
	\dfrac{\mathrm{d}^2 \boldsymbol{A}}{\mathrm{d} \boldsymbol{\rho}^2}=
	\left[
	\boldsymbol{0};
	\boldsymbol{0};
	\left[
	{\begin{array}{ccc}
			0 & \frac{\sqrt{3}}{2} & -\frac{\sqrt{3}}{2} \\
			\frac{\sqrt{3}}{2} & -\frac{\sqrt{3}}{4} & \frac{\sqrt{3}}{4} \\
			-\frac{\sqrt{3}}{2} & \frac{\sqrt{3}}{4} & \frac{\sqrt{3}}{4} 
	\end{array}} \right]
	\right]
\end{equation}
then the stress matrix is
\begin{equation}
	\boldsymbol{\omega} \dfrac{\mathrm{d}^2 \boldsymbol{A}}{\mathrm{d} \boldsymbol{\rho}^2} = b_1 \left[
	{\begin{array}{ccc}
			0 & \frac{\sqrt{3}}{2} & -\frac{\sqrt{3}}{2} \\
			\frac{\sqrt{3}}{2} & -\frac{\sqrt{3}}{4} & \frac{\sqrt{3}}{4} \\
			-\frac{\sqrt{3}}{2} & \frac{\sqrt{3}}{4} & \frac{\sqrt{3}}{4} 
	\end{array}} \right] =\boldsymbol{\omega}_s \dfrac{\mathrm{d}^2 \boldsymbol{A}}{\mathrm{d} \boldsymbol{\rho}^2}
\end{equation}
Here load will not affect the stability. The quadratic form with respect to the first-order flex is
\begin{equation}
	\boldsymbol{\rho}'^T \boldsymbol{\omega}  \dfrac{\mathrm{d}^2 \boldsymbol{A}}{\mathrm{d} \boldsymbol{\rho}^2} \boldsymbol{\rho}'= \frac{\sqrt{3}}{2} b_1 a_1^2
\end{equation}
therefore if $b_1>0$, i.e, with a positive self-stress, this vertex is prestress stable. 

It would be valuable to give some physical intuition to this result.  From figure 4 in the main text, a positive self-stress in the $z$-direction is a torque $M_3$ that induces tension on the inside of the torus, and compression on the outside: here we want to point out that such a stress distribution adds stability to an elastic system. 

Consider an initially straight elastic rod that is elastically deformed into a torus, with the ends joined: it has a tension outside and compression inside, but has no stiffness when the torus rotates about its centre-line \cite{schenk_zero_2014}. Now suppose there is a elastic stress-free torus --- this would provide a restoring force when rotating about its centre-line. The difference between the two torus shows that adding self-stress with tension inside and compression outside adds positive stiffness.

The degree-5 hole shown in figure \ref{fig: stress}(b) is prestress stable, stable under load, but not first-order rigid. The stress matrix is: 
\begin{equation}
	\boldsymbol{\omega} \dfrac{\mathrm{d}^2 \boldsymbol{A}}{\mathrm{d} \boldsymbol{\rho}^2} = \left[
	\scalemath{0.7}{\begin{array}{ccccc}
			\frac{\sqrt{2}}{4}b_1&  \frac{\sqrt{2}(-b_2+l_3)}{2}& \frac{1}{2}(-b_1-\sqrt{2}b_2+\sqrt{2}l_3)&      l_1 - \frac{\sqrt{2}}{2}l_2& \frac{1}{2}b_2-\frac{1}{2}l_3 \\
			\frac{\sqrt{2}(-b_2+l_3)}{2} &-\frac{\sqrt{2}}{2}b_1+l_3& -\frac{\sqrt{2}}{2}b_1+l_3 & \sqrt{2}l_1-l_2&  \frac{1}{2}(b_1-\sqrt{2}l_3) \\
			\frac{1}{2}(-b_1-\sqrt{2}b_2+\sqrt{2}l_3) & -\frac{\sqrt{2}}{2}b_1+l_3 & b_2& 0 & \frac{1}{2}(b_1-\sqrt{2}b_2)\\
			l_1 - \frac{\sqrt{2}}{2}l_2 & \sqrt{2}l_1-l_2 & 0 & b_2 & l_1-\frac{\sqrt{2}}{2}l_2-l_4 \\
			\frac{1}{2}b_2-\frac{1}{2}l_3 & \frac{1}{2}(b_1-\sqrt{2}l_3) & \frac{1}{2}(b_1-\sqrt{2}b_2) &  l_1-\frac{\sqrt{2}}{2}l_2-l_4 & \frac{1}{2}(-\sqrt{2}b_1+b_2+l_3)
	\end{array}} \right]
\end{equation}
The quadratic form is
\begin{equation}
	\boldsymbol{\rho}'^T \boldsymbol{\omega}  \dfrac{\mathrm{d}^2 \boldsymbol{A}}{\mathrm{d} \boldsymbol{\rho}^2} \boldsymbol{\rho}'= \dfrac{b_2-l_3}{4} a_1^2
\end{equation}
that is to say, if $b_2>l_3$, this single-hole is prestress stable. Here load will affect the stability.

\subsection{Second-order rigid but not prestress stable}
\label{section: example second order}

Statement (5) in proposition 4 provides some hints to find an example that is second-order rigid but not pre-stress stable. We need to require the dimension of self-stress and first-order flex to be at least 2. Consider a degree-$6$ hole with two-dimensional first-order flex and self-stress in figure \ref{fig: example3}(a). The coordinates are $O_1(0,0,0)$, $O_2(1,0,0)$, $O_3(1,1,0)$, $O_4(0,1,0)$, $A(-\sqrt{2}/2,-\sqrt{2}/2,0)$, $B(1,-\sqrt{2}/2,0)$, $C(1+\sqrt{2}/2,0,0)$, $D(1+\sqrt{2}/2,1,0)$, $E(0,1,1)$, $F(-\sqrt{2}/2,1+\sqrt{2}/2,0)$. The rigidity matrix is:

\begin{figure}[!tb]
	\noindent \begin{centering}
		\includegraphics[width=1\linewidth]{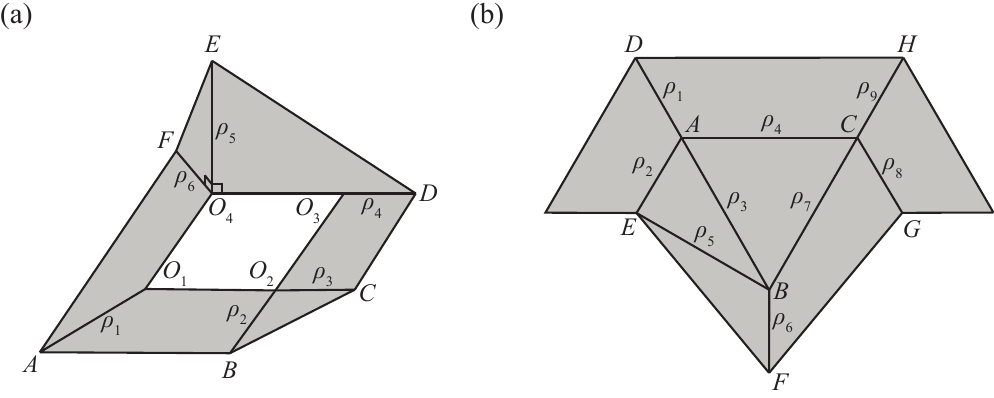}
		\par\end{centering}
	
	\caption{\label{fig: example3}(a) is a degree-6 single-hole with two-dimensional first-order flex and self-stress. (b) consists of three flat-foldable degree-4 vertices and has three-dimensional first-order flex and self-stress. (a) is examined to be prestress stable; (b) is rigid but not second-order rigid.}
\end{figure}

\begin{equation}
	\dfrac{\mathrm{d} \boldsymbol{A}}{\mathrm{d} \boldsymbol{\rho}}=
	\left[
	{\begin{array}{cccccc}
			-\frac{\sqrt{2}}{2} & 0 & 1 & 1 & 0 & -\frac{\sqrt{2}}{2}\\
			-\frac{\sqrt{2}}{2} & -1 & 0 & 0 & 0 & \frac{\sqrt{2}}{2}\\
			0 & 0 & 0 & 0 & 1 & 0 \\
			0 & 0 & 0 & 0 & 1 & 0 \\
			0 & 0 & 0 & 0 & 0 & 0 \\
			0 & -1 & 0 & -1 & 0 & \frac{\sqrt{2}}{2}
	\end{array}} \right]
\end{equation}
Since $\mathrm{rank}(\mathrm{d} \boldsymbol{A}/\mathrm{d} \boldsymbol{\rho})=4$, the first-order flex is of dimension 2, and the self-stress is of dimension 2. A basis of first-order flex is:
\begin{equation}
	[\rho']=
	\left[
	{\begin{array}{cc}
			0 &\frac{2 \sqrt{30}}{15} \\
			\frac{1}{2} & - \frac{\sqrt{15}}{10} \\
			\frac{1}{2} & \frac{\sqrt{15}}{30} \\
			0 & \frac{2 \sqrt{15}}{15} \\
			0 & 0 \\
			\frac{\sqrt{2}}{2} &\frac{\sqrt{30}}{30} \\
	\end{array}} \right]
\end{equation}
The self-stress is
\begin{equation}
	\boldsymbol{\omega}=
	\left[
	{\begin{array}{cccccc}
			0 & 0 & \frac{\sqrt{2}}{2} b_1 & -\frac{\sqrt{2}}{2} b_1 & b_2 & 0 \\
	\end{array}} \right]
\end{equation}
$b_1, b_2 \in \mathbb{R}$. The quadratic form is
\begin{equation}
	[\rho']^T \boldsymbol{\omega}_s  \dfrac{\mathrm{d}^2 \boldsymbol{A}}{\mathrm{d} \boldsymbol{\rho}^2} [\rho']= \left[
	{\begin{array}{cc}
			\frac{\sqrt{2} b_1 -2 b_2}{8} &	\frac{\sqrt{15}(\sqrt{2} b_1 -2 b_2)}{120} \\
			\frac{\sqrt{15}(\sqrt{2} b_1 -2 b_2)}{120} & \frac{\sqrt{2} b_1+30 b_2}{120}
	\end{array}} \right]
\end{equation}
From the Sylvester's criterion, $b_1 > 0$ and $0 < b_2 < \sqrt{2}b_1/2$ will stabilize this single-hole. We have tried to adjust sector angles and keep $\mathrm{rank}(\mathrm{d} \boldsymbol{A}/\mathrm{d} \boldsymbol{\rho})=4$; also we tried a rigid planar degree-6 hole which has 3 dimensional fist-order flex and self-stress, but the results are still prestress stable. Hence we conjecture that a rigid but not regular (\textit{regular} means the rigidity matrix has maximum rank) single-hole is prestress stable but not first-order rigid.

\subsection{Rigid but not second-order rigid} \label{section: example rigid}

From proposition 2 and conjecture 1 in the main text, a single-vertex or hole cannot be rigid but not second-order rigid. A possible way is to consider a rigid origami consisting several rigid-foldable inner vertices or holes, but the sector angles are carefully designed such that the intersection of their configuration space is a point, hence the whole rigid origami is rigid. 

An example is given in figure \ref{fig: example3}(b). The coordinates are $A(0,0,0)$, $B(1,-\sqrt{3},0)$, $C(2,0,0)$, $D(-1/2,\sqrt{3}/2,0)$, $E(-1/2,-\sqrt{3}/2,0)$, $F(1,-\sqrt{3}-1,0)$, $G(5/2,-\sqrt{3}/2,0)$, $H(5/2,\sqrt{3}/2,0)$. The rigidity matrix is:
\begin{equation}
	\dfrac{\mathrm{d} \boldsymbol{A}}{\mathrm{d} \boldsymbol{\rho}}=
	\left[
	{\begin{array}{ccccccccc}
			-\frac{1}{2} &	-\frac{1}{2} & \frac{1}{2} &	1  &	0&	0&	0&	0&	0 \\
			\frac{\sqrt{3}}{2}& -\frac{\sqrt{3}}{2} &	-\frac{\sqrt{3}}{2} &	0 &	0&	0&	0&	0&	0 \\
			0&	0&	0&	0&	0&	0&	0&	0&	0 \\
			0&	0&	-\frac{1}{2} &	0&	-\frac{\sqrt{3}}{2} &	0 &	\frac{1}{2} &	0&	0 \\
			0&	0&	\frac{\sqrt{3}}{2} &	0&	\frac{1}{2} &	-1 & \frac{\sqrt{3}}{2} &	0&	0 \\
			0&	0&	0&	0&	0&	0&	0&	0&	0 \\
			0&	0&	0&	-1 &	0&	0&	-\frac{1}{2} &	\frac{1}{2} &	\frac{1}{2} \\
			0&	0&	0&	0 &	0&	0&	-\frac{\sqrt{3}}{2} &	-\frac{\sqrt{3}}{2} &	\frac{\sqrt{3}}{2} \\
			0&	0&	0&	0&	0&	0&	0&	0&	0 \\
	\end{array}} \right]
\end{equation}
Since $\mathrm{rank}(\mathrm{d} \boldsymbol{A}/\mathrm{d} \boldsymbol{\rho})=6$, the first-order flex is of dimension 3, and the self-stress is also of dimension 3. A bases of first-order flex is:
\begin{equation}
	[\rho']=
	\left[
	{\begin{array}{ccc}
			- \frac{2}{\sqrt{13}} & \frac{11}{8\sqrt{13}} & \frac{\sqrt{3/2}}{8} \\
			0& \frac{\sqrt{13}}{8} & -\frac{\sqrt{3/2}}{8} \\
			- \frac{2}{\sqrt{13}}& -\frac{1}{4\sqrt{13}} & \frac{\sqrt{3/2}}{4} \\
			0 & \frac{\sqrt{13}}{8} & -\frac{\sqrt{3/2}}{8} \\
			\sqrt{\frac{3}{13}} & \frac{\sqrt{3/13}}{8} & \frac{1}{8\sqrt{2}} \\
			0& 0 & \frac{1}{\sqrt{2}} \\
			\frac{1}{\sqrt{13}} & \frac{1}{8\sqrt{13}} & \frac{3\sqrt{3/2}}{8} \\
			0 & \frac{\sqrt{13}}{8} & -\frac{\sqrt{3/2}}{8} \\
			\frac{1}{\sqrt{13}} & \frac{7}{4\sqrt{13}} & \frac{\sqrt{3/2}}{4} \\
	\end{array}} \right]
\end{equation} 
The self-stress is
\begin{equation}
	\boldsymbol{\omega}=
	\left[
	{\begin{array}{ccccccccc}
			0 & 0 & b_1 & 0 & 0 & b_2 & 0 & 0& b_3 \\
	\end{array}} \right]
\end{equation}
$b_1, b_2, b_3 \in \mathbb{R}$. The quadratic form is
\begin{equation}
	\begin{gathered}
		\left[\rho'\right]^T \boldsymbol{\omega}_s  \dfrac{\mathrm{d}^2 \boldsymbol{A}}{\mathrm{d} \boldsymbol{\rho}^2} [\rho']= \left[ \scalemath{0.8} {\begin{array}{ccc} 
				-\frac{\sqrt{3} b_2}{13}& - \frac{\sqrt{3}}{208} (26b_1+2b_2-13b_3)&	\frac{6(b_1+b_2)-3b_3}{16 \sqrt{26}} \\
				- \frac{\sqrt{3}}{208} (26b_1+2b_2-13b_3) & \frac{\sqrt{3}(117b_1-2b_2+195b_3)}{1664} & \frac{3(15b_1+2b_2+25b_3)}{128 \sqrt{26}} \\
				\frac{6(b_1+b_2)-3b_3}{16 \sqrt{26}} &	\frac{3(15b_1+2b_2+25b_3)}{128 \sqrt{26}} &	- \frac{\sqrt{3}}{256} (9b_1-26b_2+15b_3)
		\end{array}} \right] 
	\end{gathered}
\end{equation}
There is no self-stress that can stabilize the rigid origami, and there exists some first-order flex that can be extended to a second-order flex, i.e. a root for the above quadratic form for all choice of self-stress. One solution is $a_1=(9 \pm 13\sqrt{3})a_2/6$, $a_3= \sqrt{26/3}a_2$.

\subsection{Special case: when some folding angles are $\pm \pi$} \label{section: special}

Here we will discuss a special case when some of the folding angles are $\pm \pi$, which correspond to $\pm \infty$ in the normalized folding angle expression. If a rigid origami is foldable and some of the folding angles are $\pm \pi$, a flex which agrees with the consistency constraint but pointing outside $\pm \pi$ would not be valid. This ``one-side'' property brings more stiffness to a rigid origami, and the criteria on first-order rigidity, prestress stability and second-order rigidity need to be modified appropriately when some folding angles are $\pm \pi$.

\begin{figure}[!tb]
	\noindent \begin{centering}
		\includegraphics[width=1\linewidth]{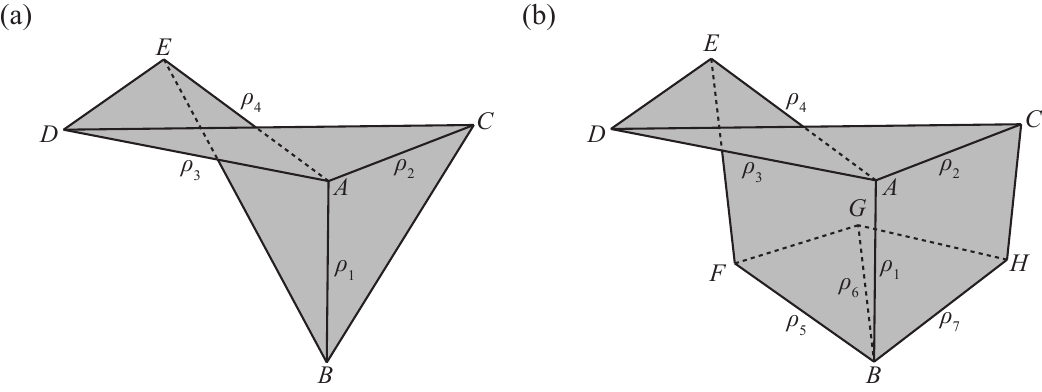}
		\par\end{centering}
	
	\caption{\label{fig: example4}(a) This is a rigid-foldable but not flat-foldable degree-4 vertex. A flex, or a first-order flex $\{\rho_1',\rho_2',\rho_3',\rho_4'\}$ obtained from the independent consistency constraint corresponds to a valid ``speed'' only when $\rho_3'>0$. If $\rho_3'<0$, the vertex will be rigid under this flex. (b) is rigid and consists of two vertices. Both of its single vertices are rigid-foldable with single folding motion. The folding angle $\rho_1$ cannot be smaller than $\pi/2$ in the top vertex $B$ and cannot be greater than $\pi/2$ in the bottom vertex $A$, hence the whole structure is rigid.}
\end{figure}

The example shown in figure \ref{fig: example4}(a) is a rigid-foldable but not flat-foldable degree-4 vertex. The sector angles are $\angle BAC=2\pi/3$, $\angle CAD=\pi/2+\arccos(-1/4)/2$, $\angle DAE=\pi/2-\arccos(-1/4)/2$, $\angle EAB=\pi/3$. A folding angle $\rho_1=\pi/2$. A global coordinate is chosen such that $B(0,0,0)$, $A(0,0,1)$, $C(\sqrt{3}/2,0,3/2)$, $D(-\sqrt{30}/10,\sqrt{30}/10,1-\sqrt{10}/5)$, $E(0,\sqrt{3}/2,1/2)$. The rigidity matrix is:

\begin{equation}
	\dfrac{\mathrm{d} \boldsymbol{A}}{\mathrm{d} \boldsymbol{\rho}}=
	\left[
	{\begin{array}{cccc}
			0 & \frac{\sqrt{3}}{2} & -\frac{\sqrt{30}}{10} & 0\\
			0 & 0 & \frac{\sqrt{30}}{10} & \frac{\sqrt{3}}{2}\\
			-1 & \frac{1}{2} & -\frac{\sqrt{10}}{5}  & -\frac{1}{2}
	\end{array}} \right]
\end{equation}

The first-order flex is:

\begin{equation}
	\boldsymbol{\rho}'=
	a_1\left[
	{\begin{array}{c}
			0 \\
			-\frac{\sqrt{2}}{3} \\
			-\frac{\sqrt{5}}{3} \\
			\frac{\sqrt{2}}{3}
	\end{array}} \right]
\end{equation}

where $a_1 \in \mathbb{R}$. Since this vertex is regular, a first-order flex is also a flex, but it corresponds to a valid ``speed'' only when $a_1<0$. If $a_1>0$, this vertex will be rigid under this flex. Such rigidity cannot be predicted by local analysis of the independent consistency constraint. 

However, this ``one-side'' property might not effect our conclusions on rigidity in some other examples. Consider a two-vertex rigid origami in figure \ref{fig: example4}(b). The sector angles are $\angle BAC=2\pi/3$, $\angle CAD=\pi/2+\arccos(-1/4)/2$, $\angle DAE=\pi/2-\arccos(-1/4)/2$, $\angle EAB=\pi/3$, $\angle ABF=\pi/2$, $\angle FBG=\pi/4$, $\angle GBH=\pi/4$, $\angle HBA=\pi/2$. A folding angle $\rho_1=\pi/2$. A global coordinate is chosen such that $A(0,0,1)$, $B(0,0,0)$, $C(\sqrt{3}/2,0,3/2)$, $D(-\sqrt{30}/10,\sqrt{30}/10,1-\sqrt{10}/5)$, $E(0,\sqrt{3}/2,1/2)$, $F(0,1,0)$, $G(\sqrt{2}/2,\sqrt{2}/2,0)$, $H(1,0,0)$. The rigidity matrix is:

\begin{equation}
	\dfrac{\mathrm{d} \boldsymbol{A}}{\mathrm{d} \boldsymbol{\rho}}=
	\left[
	{\begin{array}{ccccccc}
			0 & \frac{\sqrt{3}}{2} & -\frac{\sqrt{30}}{10} & 0 & 0 & 0 & 0\\
			0 & 0 & \frac{\sqrt{30}}{10} & \frac{\sqrt{3}}{2} & 0 & 0 & 0\\
			-1 & \frac{1}{2} & -\frac{\sqrt{10}}{5}  & -\frac{1}{2} & 0 & 0 & 0\\
			0 & 0 & 0 & 0 & 0 & \frac{\sqrt{2}}{2} & 1\\
			0 & 0 & 0 & 0 & 1 & \frac{\sqrt{2}}{2} & 0\\
			1 & 0 & 0 & 0 & 0 & 0 & 0
	\end{array}} \right]
\end{equation}

Since $\mathrm{rank}(\mathrm{d} \boldsymbol{A}/\mathrm{d} \boldsymbol{\rho})=5$, the first-order flex is of dimension 2, and the self-stress is of dimension 1. The first-order flex is:
\begin{equation}
	\boldsymbol{\rho}'=
	\left[
	{\begin{array}{c}
			0 \\
			-\frac{\sqrt{2} a_1}{3} \\
			-\frac{\sqrt{5} a_1}{3} \\
			\frac{\sqrt{2} a_1}{3} \\
			\frac{a_2}{2} \\
			-\frac{\sqrt{2}a_2}{2} \\
			\frac{a_2}{2} \\	
	\end{array}} \right]
\end{equation}
$a_1,a_2 \in \mathbb{R}$. The self-stress is
\begin{equation}
	\boldsymbol{\omega}_s=
	b_1\left[
	{\begin{array}{cccccc}
			-\sqrt{2}/4 & \sqrt{2}/4 & \sqrt{6}/4 & 0 & 0 & \sqrt{6}/4 \\
	\end{array}} \right]
\end{equation}
$b_1 \in \mathbb{R}$. The quadratic form is
\begin{equation}
	\boldsymbol{\rho}'^T \boldsymbol{\omega}_s  \dfrac{\mathrm{d}^2 \boldsymbol{A}}{\mathrm{d} \boldsymbol{\rho}^2} \boldsymbol{\rho}' = -\dfrac{\sqrt{6}b_1(9 a_1^2+10 a_2^2)}{144}
\end{equation}
It shows that the rigid origami in figure \ref{fig: example4}(b) is actually prestress stable with a negative $b_1$, even we need to require $a_2<0$.

\section{Prestress stability of a rigid planar single-vertex} \label{section: example prestress vertex}

Here we will provide the proof of proposition 2 in the main text, which states that a rigid planar single-vertex is prestress stable.

For a rigid planar degree-$n$ single-vertex, when $n \ge 4$, there must be one sector angle $\alpha_1$ greater than $\pi$ \cite{abel_rigid_2016}, which means $0< \alpha_i <\pi (2 \le i \le n)$, $0<\left(\sum_{2}^{n} \alpha_i\right)<\pi$ and all the direction vectors $\boldsymbol{x}_j (1 \le j \le n)$ lie in a semicircle. Note that $\boldsymbol{x}_1, \boldsymbol{x}_2, \cdots, \boldsymbol{x}_n$ are arranged counter-clockwise, and $\alpha_i$ is between $\boldsymbol{x}_{i-1}$ and $\boldsymbol{x}_i (2 \le i \le n)$. It turns out that for all $1 \le j<k \le n$ we have
\begin{equation}
	x_{1j}x_{2k}-x_{2j}x_{1k}=\sin \left(\sum_{j+1}^{k} \alpha_i\right)>0
\end{equation}	 

The self-stress of such single-vertex is $[0, 0, b_1], b_1 \in \mathbb{R}$. The stress matrix is 
\begin{equation}
	\boldsymbol{\omega}_s \dfrac{\mathrm{d}^2 \boldsymbol{A}}{\mathrm{d} \boldsymbol{\rho}^2} = b_1 \left[
	\scalemath{1}{\begin{array}{cccc}
			x_{11}x_{21} & x_{11}x_{22} & \cdots & x_{11}x_{2n}  \\
			x_{11}x_{22} & x_{12}x_{22} & \cdots & x_{12}x_{2n} \\
			\vdots & \vdots & \ddots & \vdots\\
			x_{11}x_{2n} & x_{12}x_{2n} & \cdots & x_{1n}x_{2n}
	\end{array}} \right]
\end{equation}

The rigidity matrix gives two equations for a first-order flex $\boldsymbol{\rho}'$ with dimension $n-2$.
\begin{equation}
	\begin{gathered}
		\sum_{1}^{n} x_{1j}\rho_j'=0 \\
		\sum_{1}^{n} x_{2j}\rho_j'=0 
	\end{gathered}
\end{equation}
hence
\begin{equation}
	\begin{gathered}
		\rho_1'=\dfrac{x_{12}\sum \limits_{3}^n x_{2j}\rho_j'-x_{22}\sum \limits_{3}^n x_{1j}\rho_j'}{x_{11}x_{22}-x_{21}x_{12}} \\
		\rho_2'=\dfrac{x_{21}\sum \limits_{3}^n x_{1j}\rho_j'-x_{11}\sum \limits_{3}^n x_{2j}\rho_j'}{x_{11}x_{22}-x_{21}x_{12}}
	\end{gathered}
\end{equation}
Write $x_{1 \times j}=x_{11}x_{2j}-x_{21}x_{1j}$, $x_{2 \times j}=x_{12}x_{2j}-x_{22}x_{1j}$, $3 \le j \le n$, the coefficient for first-order flex is actually $\boldsymbol{a}=[\rho_3'; \rho_4'; \cdots; \rho_n']$, the quadratic form is 
\begin{equation}
	\begin{gathered}
		\boldsymbol{\rho}'^T \boldsymbol{\omega}_s \dfrac{\mathrm{d}^2 \boldsymbol{A}}{\mathrm{d} \boldsymbol{\rho}^2} \boldsymbol{\rho}'= \dfrac{-b_1}{x_{1 \times 2}}\left(2\sum_{j=3}^{n} \sum_{k > j}^{n} x_{2 \times j}x_{1 \times k} \rho_j'\rho_k' + \sum_{j=3}^{n} x_{2 \times j}x_{1 \times j} \rho_j'^2 \right) \\
		=\dfrac{-b_1}{x_{1 \times 2}}\boldsymbol{a}^T\left[\begin{array}{cccc}
			x_{2 \times 3}x_{1 \times 3} & x_{2 \times 3}x_{1 \times 4} & \cdots & x_{2 \times 3}x_{1 \times n}  \\
			x_{2 \times 3}x_{1 \times 4} & x_{2 \times 4}x_{1 \times 4} & \cdots & x_{2 \times 4}x_{1 \times n} \\
			\vdots & \vdots & \ddots & \vdots\\
			x_{2 \times 3}x_{1 \times n} & x_{2 \times 4}x_{1 \times n} & \cdots & x_{2 \times n}x_{1 \times n}
		\end{array}\right]\boldsymbol{a} \\
		=\dfrac{-b_1}{\sin \alpha_2}\boldsymbol{a}^T\left[\begin{array}{cccc}
			\sin \alpha_3 \sin (\sum_2^3 \alpha_i) & \sin \alpha_3 \sin (\sum_2^4 \alpha_i) & \cdots & \sin \alpha_3 \sin (\sum_3^n \alpha_i)  \\
			\sin \alpha_3 \sin (\sum_2^4 \alpha_i) & \sin (\sum_3^4 \alpha_i) \sin (\sum_2^4 \alpha_i) & \cdots & \sin (\sum_3^4 \alpha_i) \sin (\sum_2^n \alpha_i) \\
			\vdots & \vdots & \ddots & \vdots\\
			\sin \alpha_3 \sin (\sum_3^n \alpha_i) & \sin (\sum_3^4 \alpha_i) \sin (\sum_2^n \alpha_i) & \cdots & \sin (\sum_3^n \alpha_i) \sin (\sum_2^n \alpha_i)
		\end{array}\right]\boldsymbol{a}
	\end{gathered}
\end{equation}
The $i$-th principle minor of the large matrix above is
\begin{equation}
	\left(\sin \alpha_2\right)^{i-1} \left(\prod_{3}^{i+2} \sin \alpha_j \right) \sin \left(\sum_{2}^{i+2} \alpha_j\right)>0
\end{equation}
therefore from the Sylvester's criterion, a negative self-stress $b_1<0$ will stabilize this single-vertex. Last we consider the case when $n=3$, here $0 < \alpha_i < 2\pi$ and $\alpha_1+ \alpha_2 +\alpha_3 = 2\pi$. The quadratic form from the above calculation is:
\begin{equation}
	\begin{gathered}
		\boldsymbol{\rho}'^T \boldsymbol{\omega}_s \dfrac{\mathrm{d}^2 \boldsymbol{A}}{\mathrm{d} \boldsymbol{\rho}^2} \boldsymbol{\rho}'= \dfrac{- \sin \alpha_1 \sin (\alpha_1+ \alpha_2)}{\sin \alpha_2} b_1 a_1^2
	\end{gathered}
\end{equation}
a self-stress $b_1$ such that $- b_1 \sin \alpha_1 \sin (\alpha_1+ \alpha_2)/ \sin \alpha_2 $ is positive will stabilize this single-vertex.

\section{Normalized folding angle expression} \label{section: normalized}

In order to make use of the polynomial nature behind the local rigidity concepts, in this section we provide an equivalent version of definitions of the above rigidity concepts in the normalized folding angle expression $\boldsymbol{t}=\tan (\boldsymbol{\rho}/2)$. Here the independent consistency constraint $\boldsymbol{A}(\boldsymbol{t})=\boldsymbol{0}$ becomes a polynomial system.
\begin{defn}
	A rigid origami $(\boldsymbol{t}, \lambda')$ is \textit{first-order rigid} if the only solution of $\mathrm{d} \boldsymbol{A}/\mathrm{d} \boldsymbol{t} \cdot \boldsymbol{t}'= \boldsymbol{0}$ with respect to $\boldsymbol{t}'$ is $\boldsymbol{0}$, equivalently, the rank of \textit{rigidity matrix} $\mathrm{d} \boldsymbol{A}/\mathrm{d} \boldsymbol{t}$ equals to the number of inner creases $N_c$. Otherwise this rigid origami is \textit{first-order rigid-foldable}. A non-zero $\boldsymbol{t}'$ is called a \textit{first-order flex} which forms a linear space of dimension $N_c-\mathrm{rank}(\mathrm{d} \boldsymbol{A}/\mathrm{d} \boldsymbol{t})=N_c-\mathrm{rank}(\mathrm{d} \boldsymbol{A}/\mathrm{d} \boldsymbol{\rho})$.
\end{defn}

The rigidity matrix is
\begin{equation} 
	\dfrac{\mathrm{d} \boldsymbol{A}}{\mathrm{d} \boldsymbol{t}}=\dfrac{\mathrm{d} \boldsymbol{A}}{\mathrm{d} \boldsymbol{\rho}}\dfrac{2}{1+\boldsymbol{t}^2}
\end{equation}
where $2/(1+\boldsymbol{t}^2)$ is a $N_c \times N_c$ diagonal matrix
\begin{equation} 	
	\dfrac{2}{1+\boldsymbol{t}^2}=
	\left[{\begin{array}{cccc}
			\dfrac{2}{1+t_1^2} &  &  & \\
			& \dfrac{2}{1+t_2^2} &  &  \\
			&  & \ddots &  \\
			&  &  & \dfrac{2}{1+t_{N_c}^2} \\
	\end{array}} \right] 
\end{equation}

The equilibrium is
\begin{equation} 
	\boldsymbol{\omega} \dfrac{\mathrm{d} \boldsymbol{A}}{\mathrm{d} \boldsymbol{t}} = 2\boldsymbol{l} \cdot \dfrac{1}{1+\boldsymbol{t}^2} = \dfrac{2\boldsymbol{l}}{1+\boldsymbol{t}^2}= -\dfrac{\mathrm{d} V}{\mathrm{d} \boldsymbol{t}}
\end{equation}
\begin{equation} 
	\boldsymbol{\omega}_s \dfrac{\mathrm{d} \boldsymbol{A}}{\mathrm{d} \boldsymbol{t}} = \boldsymbol{0} 
\end{equation}
when $\boldsymbol{l}$ is the work conjugate of a first-order flex $\boldsymbol{\rho}'$, $2\boldsymbol{l}/(1+\boldsymbol{t}^2)$ is the work conjugate of a first-order flex $\boldsymbol{t}'$. 
We could see that $\mathrm{d} \boldsymbol{A}/\mathrm{d}\boldsymbol{t}$ shares the same left null space with $\mathrm{d} \boldsymbol{A}/\mathrm{d}\boldsymbol{\rho}$.

\begin{thm}
	For a rigid origami $(\boldsymbol{t}, \lambda')$ with $N_v$ inner vertices, $N_h$ holes and $N_c$ inner creases, The following statements are equivalent:
	\begin{enumerate} [(1)]
		\item $(\boldsymbol{t}, \lambda')$ is first-order rigid.
		\item $(\boldsymbol{t}, \lambda')$ is statically rigid.
		\item The dimension of the collection of self-stress at $(\boldsymbol{t}, \lambda')$ is $3N_v+6N_h-N_c$.
	\end{enumerate}
\end{thm}

Then the stability condition becomes
\begin{defn} 
	A rigid origami $(\boldsymbol{t},\lambda')$ with $N_v$ inner vertices, $N_c$ inner creases and $N_h$ holes is \textit{prestress stable} if there is a positive-definite matrix $B$ with size ($3N_v+6N_h) \times (3N_v+6N_h$), and a vector $\boldsymbol{\omega}_s \in \mathbb{R}^{3N_v+6N_h}$ such that
	\begin{equation} 
		\boldsymbol{\omega}_s \dfrac{\mathrm{d} \boldsymbol{A}}{\mathrm{d} \boldsymbol{t}}=\boldsymbol{0}
	\end{equation}
	and
	\begin{equation} 
		K=\dfrac{\mathrm{d} \boldsymbol{A}}{\mathrm{d} \boldsymbol{t}}^T  B \dfrac{\mathrm{d} \boldsymbol{A}}{\mathrm{d} \boldsymbol{t}}+\boldsymbol{\omega}_s \dfrac{{\mathrm{d}^2 \boldsymbol{A}}}{\mathrm{d} \boldsymbol{t}^2}
	\end{equation}
	is positive-definite.  
\end{defn}

The Hessian of independent constraint in the normalized folding angle expression is
\begin{equation}
	\dfrac{\mathrm{d}^2 A_i}{ \mathrm{d} t_k \mathrm{d} t_j}=\dfrac{\mathrm{d}^2 A_i}{ \mathrm{d} \rho_k \mathrm{d} \rho_j}\dfrac{4}{(1+t_k^2)(1+t_j^2)}-\dfrac{\mathrm{d} A_i}{ \mathrm{d} \rho_j}\dfrac{4\delta_{kj}t_j}{(1+t_j^2)^2}
\end{equation}
where $\delta$ is the Kronecker delta, and therefore
\begin{equation}
	\boldsymbol{\omega}_s \dfrac{{\mathrm{d}^2 \boldsymbol{A}}}{\mathrm{d} \boldsymbol{t}^2}=\omega_i \dfrac{\mathrm{d}^2 A_i}{ \mathrm{d} t_k \mathrm{d} t_j}= \omega_i \dfrac{\mathrm{d}^2 A_i}{ \mathrm{d} \rho_k \mathrm{d} \rho_j}\dfrac{4}{(1+t_k^2)(1+t_j^2)}= \dfrac{2}{1+\boldsymbol{t}^2} \boldsymbol{\omega}_s \dfrac{{\mathrm{d}^2 \boldsymbol{A}}}{\mathrm{d} \boldsymbol{\rho}^2} \dfrac{2}{1+\boldsymbol{t}^2}
\end{equation}

Consider reducing the calculation:
\begin{prop} 
	\begin{enumerate} [(1)]
		\item The matrix $[(\mathrm{d} \boldsymbol{A}/\mathrm{d} \boldsymbol{t})^T \cdot B \cdot \mathrm{d} \boldsymbol{A}/\mathrm{d} \boldsymbol{t}]$ is positive semi-definite.  Its quadratic and linear nullspace are equal to the collection of all first-order flex that form the nullspace of $\mathrm{d} \boldsymbol{A}/\mathrm{d} \boldsymbol{t}$.
		\item A rigid origami $(\boldsymbol{t},\lambda')$ is prestress stable if and only if there exists a self-stress $\boldsymbol{\omega}_s \in \mathbb{R}^{3N_v+6N_h}$ such that $\boldsymbol{\omega}_s \cdot \mathrm{d}^2 \boldsymbol{A}/\mathrm{d} \boldsymbol{t}^2$ is positive-definite when restricted to the nullspace of $\mathrm{d} \boldsymbol{A}/\mathrm{d} \boldsymbol{t}$. 
		\item A rigid origami $(\boldsymbol{t},\lambda')$ is prestress stable if and only if there exists a self-stress $\boldsymbol{\omega}_s \in \mathbb{R}^{3N_v+6N_h}$ such that all the eigenvalues of $[t']^T \cdot \boldsymbol{\omega}_s \cdot \mathrm{d}^2 \boldsymbol{A}/\mathrm{d} \boldsymbol{t}^2 \cdot [t']$ are positive, where $[t']$ is the collection of a basis of first-order flex.
		\begin{equation} 
			[t']= \left[\begin{array}{cccc}
				\boldsymbol{t}_1' & \boldsymbol{t}_2' & \cdots & \boldsymbol{t}_{N_c-\mathrm{rank}(\mathrm{d} \boldsymbol{A}/\mathrm{d} \boldsymbol{t})}'
			\end{array} \right]
		\end{equation}
	\end{enumerate}
\end{prop}

The second-order rigidity condition is hence
\begin{defn} 
	For a rigid origami $(\boldsymbol{t}',\lambda')$ with $N_v$ inner vertices, $N_c$ inner creases and $N_h$ holes, a \textit{second-order flex} $(\boldsymbol{t}',\boldsymbol{t}'') \in (\mathbb{R}^{N_c},\mathbb{R}^{N_c})$ satisfies 
	\begin{equation} \label{eq: second-order t}
		\begin{gathered}
			\dfrac{\mathrm{d} \boldsymbol{A}}{\mathrm{d} \boldsymbol{t}} \boldsymbol{t}'=\boldsymbol{0} \\
			\boldsymbol{t}'^T  \dfrac{\mathrm{d}^2 \boldsymbol{A}}{\mathrm{d} \boldsymbol{t}^2}  \boldsymbol{t}' + \dfrac{\mathrm{d} \boldsymbol{A}}{\mathrm{d} \boldsymbol{t}}  \boldsymbol{t}''=\boldsymbol{0}
		\end{gathered}
	\end{equation}
	A second-order flex with $\boldsymbol{t}'=\boldsymbol{0}$ is called \textit{trivial}, otherwise \textit{non-trivial}. If there is only trivial second-order flex, this rigid origami is \textit{second-order rigid}, otherwise \textit{second-order rigid-foldable}.
\end{defn}

\begin{prop} 
	\begin{enumerate}[(1)]
		\item A first-order flex $\boldsymbol{t}'$ can be extended to a second-order flex $\boldsymbol{t}''$ if and only if for all self-stress $\boldsymbol{\omega}_s$,  $\boldsymbol{t}'^T  [\boldsymbol{\omega}_s \cdot \mathrm{d}^2 \boldsymbol{A}/\mathrm{d} \boldsymbol{t}^2]\, \boldsymbol{t}'= 0$.
		\item A rigid origami is second-order rigid if and only if for any first-order flex $\boldsymbol{t}'$ there is a self-stress $\boldsymbol{\omega}_s(\boldsymbol{t}')$ s.t.,  $\boldsymbol{t}'^T [\boldsymbol{\omega}_s \cdot \mathrm{d}^2 \boldsymbol{A}/\mathrm{d} \boldsymbol{t}^2]\, \boldsymbol{t}' > 0$.
		\item A rigid origami is second-order rigid if and only if for any first-order flex $\boldsymbol{t}'$ there is a self-stress $\boldsymbol{\omega}_s(\boldsymbol{t}')$ s.t.,  $\boldsymbol{t}'^T [\boldsymbol{\omega}_s \cdot \mathrm{d}^2 \boldsymbol{A}/\mathrm{d} \boldsymbol{t}^2]\, \boldsymbol{t}' > 0$.
		\item A rigid origami is second-order rigid if and only if the intersection of quadratic nullspace of all $[t']^T \cdot \boldsymbol{\omega}_i \cdot \mathrm{d}^2 \boldsymbol{A}/\mathrm{d} \boldsymbol{t}^2 \cdot [t']$ is $\boldsymbol{0}$. Here $\{\boldsymbol{\omega}_i\}$ is a base of self-stress $(1 \le i \le 3N_v+6N_h-\mathrm{rank}(\mathrm{d} \boldsymbol{A}/\mathrm{d} \boldsymbol{t}))$.
	\end{enumerate}
\end{prop}

We claim that the above first-order or static rigidity, prestress stability and second-order rigidity are equivalent for the normalized folding angle expression $\boldsymbol{t}=\tan (\boldsymbol{\rho}/2)$ and the folding angle expression $\boldsymbol{\rho}$. Such equivalence is presented below. 

\begin{prop}
	Given a rigid origami $(\boldsymbol{\rho},\lambda')$, which can also be written as $(\boldsymbol{t},\lambda')$, the following conclusions hold.
	\begin{enumerate} [(1)]
		\item $\boldsymbol{\rho}'$ is a first-order flex for $(\boldsymbol{\rho},\lambda')$ if and only if $\boldsymbol{t}'=\tan (\boldsymbol{\rho}'/2)$ is a first-order flex for $(\boldsymbol{t},\lambda')$.
		\item $\boldsymbol{\omega}$ is a stress for $(\boldsymbol{\rho},\lambda')$ under load $\boldsymbol{l}$ if and only if $\boldsymbol{\omega}$ is a stress for $(\boldsymbol{t},\lambda')$ under load $2\boldsymbol{l}/(1+\boldsymbol{t}^2)$.
		\item $\boldsymbol{\omega}_s$ is a self-stress for $(\boldsymbol{\rho},\lambda')$ if and only if $\boldsymbol{\omega}_s$ is a stress for $(\boldsymbol{t},\lambda')$.
		\item $\boldsymbol{\omega}$ stabilizes $(\boldsymbol{\rho},\lambda')$ if and only if $\boldsymbol{\omega}$ stabilize $(\boldsymbol{t},\lambda')$.
		\item $\boldsymbol{\omega}_s$ stabilizes $(\boldsymbol{\rho},\lambda')$ if and only if $\boldsymbol{\omega}_s$ stabilize $(\boldsymbol{t},\lambda')$.
		\item $\boldsymbol{\omega}_s$ blocks a first-order flex $\boldsymbol{\rho}'$ for $(\boldsymbol{\rho},\lambda')$ if and only if $\boldsymbol{\omega}_s$ blocks a first-order flex $\boldsymbol{t}'$ for $(\boldsymbol{t},\lambda')$.
	\end{enumerate}
\end{prop}  

\bibliographystyle{unsrt}